\input amstex
\input xy
\xyoption{all}
\SelectTips {cm}{}

\documentstyle{amsppt}
\pagewidth{35pc}

\topmatter
\title
Manifolds, Structures Categorically
\endtitle
\author
G.V. Kondratiev
\endauthor
\address
Department of Mathematics, Ottawa University
\endaddress
\email
gennadii\@\, hotmail.com
\endemail
\subjclass
CT Category Theory
\endsubjclass
\abstract
A notion of general manifolds is introduced. It covers all usual manifolds in mathematics. Essentially, it is a way
how to get a bigger 'fibration' over a site which locally coincides with a given one. An enrichment with generalized elements
is regarded which allows to see hom-sets of a given category as (almost) objects and to transfer some technics from objects onto
hom-sets. Lifting problem for a group action and actions of group objects are also included.
\endabstract

\endtopmatter

\document

\head {\bf 1. Fibrations and cofibrations} \endhead

(Co)fibrations play role of structures over objects in a given category which can be transported along morphisms.
Transport system is called (co)cartesian morphisms. The situation is very similar to fibrations with connection 
as in Differential Geometry.

\vskip 0.2cm
{\bf Definition 1.1.} For a functor \ $\vcenter{\xymatrix{\bold {E} \ar[d]_-{p} \\ \bold {B}}}$ 
\item{$\bullet $} morphism $f:B'\to p(E)$ has a {\bf cartesian lifting} \ $\tilde f:E'\to E\in Ar\, \bold {E}$ \ if $\forall f':E''\to E$ 
such that $p(f')$ factors through $f$ (i.e., $\exists g:p(E'')\to B'$ such that $p(f')=f\circ g$) \ $f'$ itself uniquely factors 
through $\tilde f$ over the base factorization (i.e., $\exists !\, \tilde g:E''\to E'$ such that $f'=\tilde f\circ \tilde g$ and $p(\tilde g)=g$) \cite{Jac} 
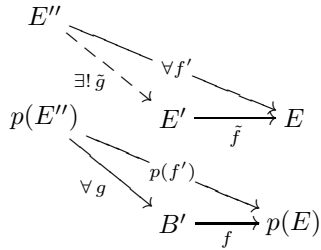
$$\xymatrix{E'' \ar[rrd]|-{\ \forall f'} \ar@{-->}[rd]_-{\exists !\, \tilde g} &&\\
p(E'') \ar[rrd]|-{\ p(f')\ } \ar[rd]_-{\forall \, g} & E' \ar[r]_-{\tilde f} & E \\
 & B' \ar[r]_-{f} & p(E)
}$$
\item{$\bullet $} morphism $f:p(E)\to B'$ has a {\bf cocartesian lifting} \ $\tilde f:E\to E'\in Ar\, \bold {E}$ \ if $\forall f':E\to E''$ 
such that $p(f')$ factors through $f$ (i.e., $\exists g:B'\to p(E'')$ such that $p(f')=g\circ f$) \ $f'$ itself uniquely factors 
through $\tilde f$ over the base factorization (i.e., $\exists !\, \tilde g:E'\to E''$ such that $f'=\tilde g\circ \tilde f$ and $p(\tilde g)=g$) \cite{Jac}
$$\xymatrix{&& E'' \\
E \ar[r]_-{\tilde f} \ar[rru]|-{\ \forall f'\ } & E' \ar@{-->}[ru]_-{\exists !\, \tilde g} & p(E'') \\
p(E) \ar[r]_-{f} \ar[rru]|-{\ \ p(f')\ \, } & B' \ar[ru]_-{\forall \, g} &
}$$  \vskip -0.2cm \hfill $\square $
\vskip 0.2cm
{\bf Remark.} (Co)cartesian morphisms if exist are unique up to vertical isomorphism ($v:E\to E'\in Ar\, \bold{E}$ is {\bf vertical} if 
$p(v)=1_B$ for some $B\in Ob\, \bold{B}$).

\vskip 0.2cm
{\bf Definition 2.1.2.} A functor \ $\vcenter{\xymatrix{\bold {E} \ar[d]_-{p} \\ \bold {B}}}$ \ is called \cite{Jac, Str}
\item{$\bullet $} {\bf fibration} \ if for each \ $f:B'\to p(E)\in Ar\, \bold {B}$ \ there exists cartesian lifting \ $\tilde f:E'\to E\in Ar\, \bold {E}$
\item{$\bullet $} {\bf cofibration} \ if for each \ $f:p(E)\to B'\in Ar\, \bold {B}$ \ there exists cocartesian lifting \ $\tilde f:E\to E'\in Ar\, \bold {E}$
\item{$\bullet $} {\bf bifibration} \ if it is both fibration and cofibration    \hfill $\square $

\vskip 0.2cm
Subcategory $\bold {E}_B:=p^{-1}(B,1_B)\hookrightarrow \bold {E}$ is called {\bf fiber} over $B$. $\bold{E}_B$ consists of all vertical morphisms over $B$.

\vskip 0.4cm
\centerline{\bf Examples}

\item{1.} For a category $\bold{C}$ with pullbacks {\bf codomain fibration} is \ 
$\vcenter{\xymatrix{\bold{\ \ C^{\to }} \ar[d]_-{\text{cod}} \\ \bold{C}}}$, \ \ 
$\text{cod}\Bigl(\vcenter{\xymatrix{\bullet \ar[r]^-{\Phi } \ar[d] & \bullet \ar[d] \\ \bullet \ar[r]_-{f} & \bullet }}\Bigr)=f$, \ \ 
cartesian lifting is a pullback square
\item{2.} For a category $\bold{C'}$ with pushouts {\bf domain cofibration} is \ 
$\vcenter{\xymatrix{\bold{\ \ \ C^{'\to }} \ar[d]_-{\text{dom}} \\ \bold{C'}}}$, \ \ 
$\text{dom}\Bigl(\vcenter{\xymatrix{\bullet \ar[r]^-{\Phi } \ar[d] & \bullet \ar[d] \\ \bullet \ar[r]_-{f} & \bullet }}\Bigr)=\Phi$, \ \
cocartesian lifting is a pushout square
\item{3.} \cite{Kom} \ Denote $\bold{Rng-}\bold{Mod}$, category of left modules with variable rings. \newline
$Ob\, (\bold{Rng-}\bold{Mod})$ are pairs $(R,M)$ where: $R$ is a ring, $M$ is an $R$-module.
$Ar\, (\bold{Rng-}\bold{Mod})$ are pairs $(\varphi :R_1\to R_2,f:M_1\to M_2)$ such that
$\forall r\in R_1,m\in M_1 \ \ f(r\cdot m)=\varphi (r)\cdot f(m)$.
Then 'projection on first component' $\bold{Rng-}\bold{Mod}@>p_1>>\bold{Rng}$ is a bifibration.
If $\varphi :R_1\to R_2\in Ar\, \bold{Rng}$ then
{\bf direct image} of $(R_1,M_1)$ over $\varphi $ is $(R_2,R_2\otimes _{R_1}M_1)$.
Cocartesian lifting of $\varphi $ is $(\varphi ,\bar \varphi )$ with $\bar \varphi :M_1\to R_2\otimes _{R_1}M_1:m\mapsto 1\otimes m$.

$$\xymatrix{ && (R_3,M_3)\\
(R_1,M_1) \ar[urr]^-{(\psi \circ \varphi ,f)} \ar[r]_-{(\varphi ,\bar \varphi )} & (R_2,R_2\otimes _{R_1}M_1) \ar@{-->}[ur]_-{(\psi ,\psi \cdot f)}\\
 && R_3\\
R_1 \ar[r]_-{\varphi } \ar[urr]^-{\psi \circ \varphi } & R_2 \ar[ur]_-{\psi } &
}$$

where: $\psi \cdot f:=\mu \circ (\psi \otimes f):$
$R_2\otimes _{R_1}M_1@>{\psi \otimes f}>>R_3\otimes _{R_1}M_3@>\mu >>M_3$, \
$\mu $ is module multiplication, \

$\psi (r_1\cdot r_2)\cdot f(m)=(\psi \cdot f)(r_1\cdot r_2\otimes m)=\psi (r_1)\cdot (\psi \cdot f)(r_2\otimes m)=\psi (r_1)\cdot (\psi (r_2)\cdot f(m))$.

{\bf Inverse image} of $(R_1, M_1)$ over $\alpha :R\to R_1$ is $(R,M_1)$ with action $(r,m)\mapsto \alpha (r)\cdot m$.

Cartesian lifting of $\alpha $ is $(\alpha , 1_{M_1})$.

$$\xymatrix{(R',M') \ar[drr]^-{(\alpha \circ \beta ,f)} \ar@{-->}[dr]_-{(\beta ,f)} &&\\
 & (R,M_1) \ar[r]_-{(\alpha ,1_{M_1})} & (R_1,M_1)\\
R' \ar[drr]^-{\alpha \circ \beta } \ar[dr]_-{\beta } &&\\
 & R \ar[r]_-{\alpha } & R_1
}$$
\vskip 0.35cm
\item{4.} (Differential Equations)\vskip -0.00cm 
$\xymatrix{\bold{Set} & \bold{Top} \ar[l]_-{p_1} & \bold{Diff} \ar[l]_-{p_2} & \bold{Diff}^{\twoheadrightarrow } \ar[l]_-{p_3} & \bold{Jet}^{\infty }(\bold{Diff}^{\twoheadrightarrow }) \ar[l]_-{p_4} & \bold{Sub}(\bold{Jet}^{\infty }(\bold{Diff}^{\twoheadrightarrow })) \ar[l]_-{p_5}}$
\vskip 0.1cm
where $\bold{Diff}^{\twoheadrightarrow }$ is a full subcategory of $\bold{Diff}^{\to }$ consisting of surjective submersions 
(surjective maps with surjective differential at each point), $\bold{Jet}^{\infty }(\bold{Diff}^{\twoheadrightarrow })$ is the 
corresponding category of $\infty $-jet-bundles, $\bold{Sub}(\bold{Jet}^{\infty }(\bold{Diff}^{\twoheadrightarrow }))$ is the category
of subobjects of $\infty $-jet-bundles (differential relations),
\vskip -0.3cm
$p_1:\bold{Top}\to \bold{Set}$ is a fibration (\, $\forall $ arrow $S'@>u>>p_1(T)$ \, $\exists $ \, Cartesian 'completion' $\vcenter{\xymatrix{T' \ar[r]^-{v} & T \\
S' \ar[r]^-{u} & p_1(T)}}$, $T'$ has initial topology w.r.t. $u$),
\vskip 0.0cm
$p_2:\bold{Diff}\to \bold{Top}$ is not a fibration if differentiable manifolds are regarded in usual sense (as locally Eucledian), but
it is a fibration if differentiable manifolds are topological spaces endowed with a subsheaf of continuous functions closed under 
smooth operations,
\vskip 0.0cm
$p_3:\bold{Diff}^{\twoheadrightarrow }\to \bold{Diff}$ is a codomain fibration since pullback of a surjective submersion is a surjective 
submersion,
\vskip 0.0cm
$p_4:\bold{Jet}^{\infty }(\bold{Diff}^{\twoheadrightarrow })\to \bold{Diff}^{\twoheadrightarrow }$ is not a fibration (if we admit arbitrary fibre bunble arrows between objects in $\bold{Jet}^{\infty }(\bold{Diff}^{\twoheadrightarrow })$), but it is a structure
over $\bold{Diff}^{\twoheadrightarrow }$ (see 2.1.1),
\vskip 0.0cm
$p_5:\bold{Sub}(\bold{Jet}^{\infty }(\bold{Diff}^{\twoheadrightarrow }))\to \bold{Jet}^{\infty }(\bold{Diff}^{\twoheadrightarrow })$ is
a 'subobject' fibration.

\vskip 0.35cm
\proclaim{\bf Lemma 1.1} Every functor $F:\bold{C}\to \bold{B}$ factors through a {\bf free fibration} $\xymatrix{\bold{C} \ar[r]^-{i} \ar[dr]_-{F} & 1/F \ar[d]^-{\text{dom}} \\
 & \bold{B}
}$
where $i:\bold{C}\to 1/F:\cases C\mapsto 1_{F(C)} & C\in Ob\, \bold{C}\\
(f:C\to C')\mapsto (F(f),f) & f\in Ar\, \bold{C}\endcases $   \vskip 0.05cm \hfill $\square $
\endproclaim

\vskip 0.0cm
Taking (co)cartesian morphism over $f:B'\to B\in Ar\, \bold{B}$ depends on the choice of an object $E\in Ob\, \bold{E}_B$ 
(respectively, $E'\in Ob\, \bold{E}_{B'}$) and the choice of an arrow $\tilde f_E:E'\to E$ (respectively, $\tilde f_{E'}:E'\to E$) in 
the isomorphism class.

\vskip 0.2cm
\proclaim{\bf Lemma 1.2} \item{$\bullet $} For a (co)fibration $\vcenter{\xymatrix{\bold{E} \ar[d]_-{p} \\ \bold{B}}}$, an arrow $f:B'\to B\in Ar\, \bold{B}$, and a choice $\tilde f_E:E'\to E$
\ $\forall E\in Ob\, \bold{E}_B$ (respectively, $\tilde f_E':E'\to E$ \ $\forall E'\in Ob\, \bold{E}_{B'}$) there is a functor \ 
$\bold{Cart}_f:\bold{E}_B\to \bold{E}_{B'}:$
$$\cases E\mapsto d(\tilde f_E) & E\in Ob\, \bold{E}_B\\
(g:E_1\to E)\mapsto \bold{Cart}_f(g) \ \ \text{ s.t. }  \vcenter{\xymatrix{E'_1 \ar[r]^-{\tilde f_{E_1}} \ar@{-->}[d]_-{\exists !\, \bold{Cart}_f(g)} & E_1 \ar[d]^-{g} \\
E' \ar[r]_-{\tilde f_E} & E}}        &   g\in Ar\, \bold{E}_B
\endcases $$ 
respectively, \ $\bold{Cart}_f:\bold{E}_{B'}\to \bold{E}_{B}:$
$$
\cases E'\mapsto c(\tilde f_{E'}) & E'\in Ob\, \bold{E}_{B'}\\
(g:E'_1\to E')\mapsto \bold{Cart}_f(g) \ \ \text{ s.t. }  \vcenter{\xymatrix{E'_1 \ar[r]^-{\tilde f_{E'_1}} \ar[d]_-{g} & E_1 \ar@{-->}[d]^-{\exists !\, \bold{Cart}_f(g)} \\
E' \ar[r]_-{\tilde f_{E'}} & E}}        &   g\in Ar\, \bold{E}_{B'}
\endcases $$
\item{$\bullet $} For any two such choices of (co)cartesian morphisms $\tilde f_E, \ \forall E\in \bold{E}_B$ and $\Tilde {\Tilde f}_E, \ \forall E\in \bold{E}_B$
the corresponding functors $\bold{Cart}_f$ and $\widetilde {\bold{Cart}}_f$ are isomorphic.
\endproclaim
\demo\nofrills{Proof}\ \,  is straightforward. \ See \cite{Jac}. \hfill $\square $
\enddemo

Another equivalent description of (co)fibrations is via (co)contravariant pseudofunctors \linebreak
$\bold{B}\to \bold{CAT}$.

\vskip 0.2cm
\proclaim{\bf Proposition 1.1}
\vskip 0.15cm
\item{$\bullet $} For each choice of (co)cartesian liftings for all arrows $f:B'\to p(E)$ (or respectively, $f:p(E)\to B'$)
in the base category of a (co)fibration $\vcenter{\xymatrix{\bold{E} \ar[d]_-{p} \\ \bold{B}}}$ there is a corresponding pseudofunctor
$F_p:\bold{B}\to \bold{CAT}:\cases B\mapsto \bold{E}_B & B\in Ob\, \bold{B}\\
f\mapsto \bold{Cart}_f & f\in Ar\, \bold{B}
\endcases$
\vskip 0.1cm
\item{$\bullet $} Conversely, for a given (co)contravariant pseudofunctor $F:\bold{B}\to \bold{CAT}$ there is a (co)fibration
$\vcenter{\xymatrix{\bold{E}_F \ar[d]_-{p_F} \\ \bold{B}}}$ \ ({\bf Grothendieck construction}).
$$\cases Ob\, \bold{E}_F \text{ are pairs }
\dbinom EB,  & \ \, E\in Ob\, F(B), \ \, B\in Ob\, \bold{B}\\
Ar\, \bold{E}_F \text{ are pairs } \dbinom hf\in \bold{E}_F\Bigl (\dbinom EB,\dbinom {E'}{B'}\Bigr ) & \ \, h\in F(B)(E,F(f)(E')), \ \, f\in \bold{B}(B,B')
\endcases $$
$1_{\tbinom EB}:=\dbinom {E@>\sim >>F(1_B)(E)}{1_B}$,
\vskip 0.1cm
for $\dbinom uf:\dbinom EB\to \dbinom {E'}{B'})$ and $\dbinom vg:\dbinom {E'}{B'})\to \dbinom {E''}{B''}$ \ the composite
\vskip 0.1cm
$\dbinom vg\circ \dbinom uf:=\dbinom w{g\circ f}$
where $w:E@>u>>F(f)(E')@>F(f)(v)>>F(f)(F(g)(E''))@>\sim >>F(g\circ f)(E'')$
\vskip 0.2cm
or respectively,
$$\cases Ob\, \bold{E}_F \text{ are pairs }
\dbinom EB,  & \ \, E\in Ob\, F(B), \ \, B\in Ob\, \bold{B}\\
Ar\, \bold{E}_F \text{ are pairs } \dbinom hf\in \bold{E}_F((E,B),(E',B')) & \ \, h\in F(B')(F(f)(E),E'), \ \, f\in \bold{B}(B,B')
\endcases $$

$1_{\tbinom EB}:=\dbinom {F(1_B)(E)@>\sim >>E}{1_B}$,
\vskip 0.1cm
for $\dbinom uf:\dbinom EB\to \dbinom {E'}{B'}$ and $\dbinom vg:\dbinom {E'}{B'}\to \dbinom {E''}{B''}$ \ the composite
\vskip 0.1cm
$\dbinom vg\circ \dbinom uf:=\dbinom w{g\circ f}$
where $w:F(g\circ f)(E)@>\sim >>F(g)(F(f)(E))@>F(g)(u)>>F(g)(E')@>v>>E''$,
\vskip 0.2cm
$p_F$ is the projection onto the bottom component in both cases.

\vskip 0.1cm
\item{$\bullet $} Two processes above are (weakly) inverse to each other.
\endproclaim
\demo\nofrills{Proof} \ \ is straightforward. \ See \cite{Jac}.

(It is essential for Grothendieck construction that every (co)fibration $\vcenter{\xymatrix{\bold{E} \ar[d]_-{p} \\ \bold{B}}}$ is equivalent to (co)fibration $\vcenter{\xymatrix{p \ar[d]_-{p_2} \\ \bold{B}}}$ where $\cases Ob\, p:=\Bigl \{\dbinom E{p(E)}\ \Bigl |\Bigr . \ E\in Ob\, \bold{E}\Bigr \} & \\
Ar\, p:=\Bigl \{\dbinom f{p(f)}\ \Bigl | \Bigr .\  f\in Ar\, \bold{E}\Bigr \} & \endcases $\hskip -0.35cm, \, $p_2$ is the  
projection onto the bottom\linebreak
\vskip 0.1cm 
\hskip -0.35cm component, \, and that every morphism in $\bold{E}$ factors through (co)cartesian one)    \hfill $\square $
\enddemo

\vskip 0.3cm
\head {\bf 2. Almost structures} \endhead

{\bf Definition 2.1.} {\bf Structure of type $\bold{E}$ on {\rm (}objects of{\rm )} category $\bold{B}$} is a functor \hskip 0.35cm$\vcenter{\xymatrix{\bold{E} \ar[d]_-{p} \\ \bold{B}}}$
\vskip -0.2cm
which is
\vskip 0.1cm
\item{$\bullet $} faithful
\vskip 0.1cm
\item{$\bullet $} admits lifting of iso's of type $f:B'@>\sim >> p(E)$ (or, the same, $f:p(E)@>\sim >> B'$)
\vskip 0.1cm
\item{$\bullet $} each fiber $\bold{E}_B$ is skeletal     \hfill $\square $
\vskip 0.2cm
\proclaim{\bf Lemma 2.1} Let $\vcenter{\xymatrix{\bold{E} \ar[d]_-{p} \\ \bold{B}}}$ be a structure on \, $\bold{B}$, \, $E',E''\in Ob\, (\bold{E}_B)$ \, for some \, $B\in Ob\, (\bold{B})$.
\vskip 0.0cmThe following are equivalent
\vskip 0.1cm
\item{a)} $E'=E''$
\vskip 0.1cm
\item{b)} $\forall E\in Ob\, \bold{E}$ \, $p(\bold{E}(E,E'))=p(\bold{E}(E,E''))$
\vskip 0.1cm
\item{c)} $\forall E\in Ob\, \bold{E}$ \, $p(\bold{E}(E',E))=p(\bold{E}(E'',E))$
\vskip 0.1cm
\item{d)} $\forall E\in Ob\, \bold{E}_B$ \, $\bold{E}_B(E,E')=\emptyset $ \, iff \, $\bold{E}_B(E,E'')=\emptyset $
\vskip 0.1cm
\item{e)} $\forall E\in Ob\, \bold{E}_B$ \, $\bold{E}_B(E',E)=\emptyset $ \, iff \, $\bold{E}_B(E'',E)=\emptyset $
\endproclaim
\demo{Proof} $a)\Rightarrow b),c),d),e)$ is obvious.\newline
$b)\Rightarrow a)$ Take $E=E'$ and $E=E''$ then $\exists f:E'\to E''$ and $\exists g:E''\to E'$ such that $p(f)=1_B=p(g)$. So, $p(f\circ g)=p(g\circ f)=1_B$. 
By faithfulness of $p$ and skeletal condition on $\bold{E}_B$, $f$ and $g$ are trivial iso's.\newline
$c)\Rightarrow a)$ is the same as $b)\Rightarrow a)$.\newline 
$d)\Rightarrow a)$ Take $E=E'$ and $E=E''$ then $\exists f:E'\to E''$, $f$ is vertical, and $\exists g:E''\to E'$, $g$ is vertical. 
So, $f\circ g=1_{E''}, \, g\circ f=1_{E'}\Rightarrow E'\simeq E''\Rightarrow E'=E''$.\newline
$e)\Rightarrow a)$ is the same as $d)\Rightarrow a)$.   \hfill $\square $
\enddemo
\vskip -0.2cm
\proclaim{\bf Proposition 2.1} For structure $\vcenter{\xymatrix{\bold{E} \ar[d]_-{p} \\
\bold{B}}}$ on $\bold{B}$ each fiber $\bold{E}_B$ is a poset category with (vertical) cartesian morphisms, just identities.
\endproclaim
\demo{Proof} $\bold{E}_B$ is a preorder since $\exists $ at most one morphism between objects ($p$ is faithfull). 
$\bold{E}_B$ is a partial order since all isomorphic objects are the same (skeletal condition). 
Cartesian lift of iso is iso, so cartesian lift $\widetilde {1_B}$ is an identity (since $\bold{E}_B$ is a partial order). \hfill $\square $
\enddemo

\vskip 0.2cm
\proclaim{\bf Proposition 2.2}
\vskip 0.2cm
\item{$\bullet $} Pullback of a fibration is a fibration,
\item{$\bullet $} Pullback of a structure (on $\bold{B}$) is a structure (on $\bold{A}$) \hskip 0.35cm$\vcenter{\xymatrix{F^*\bold{E} \ar[r]^-{\bar F} \ar[d]_-{p'} & \bold{E} \ar[d]^-{p} \\
\bold{A} \ar[r]_-{F} & \bold{B}}}$
\endproclaim
\demo{Proof} 
\item{$\bullet $} Pullbacks in $\bold{CAT}$ exist and they are certain subcategories of direct products. Cartesian morphisms are 
preserved under pullbacks which is seen from the following diagram 
\vskip 0.3cm
$\xymatrix{(A'',E'') \ar[drr]^-{(u\circ w,\alpha )} \ar[dr]_-{(w,\beta )} & & \\
 & (A',E') \ar[r]_-{(u,v)} & (A,E) \\
 A'' \ar[dr]^-{\forall w} & & \\
 & A' \ar[r]^-{u} & A}$
\hskip 1.5cm$\xymatrix{E'' \ar[drr]^-{\alpha } \ar[dr]_-{\beta } & & \\
 & E' \ar[r]_-{v} & E \\
 FA'' \ar[dr]^-{Fw} & & \\
 & FA' \ar[r]^-{Fu} & FA}$
 \vskip 0.0cm
 (where: $p(E)=F(A), \ p(\alpha )=F(u\circ w)$, $v$ is cartesian over $F(u)$).
\vskip 0.0cm
\item{$\bullet $} If $\vcenter{\xymatrix{\bold{E} \ar[d]_-{p} \\
\bold{B}}}$ admits lifting of iso's of type $B'@>\sim >>p(E)$ then $\vcenter{\xymatrix{F^*\bold{E} \ar[d]_-{p'} \\
\bold{A}}}$ admits lifting of iso's of the same type (obvious).
If $p$ is faithful then $p'$ is faithful (assume, $p'(u,v)=p'(u,v_1)=u$ then $F(u)=p(v)=p(v_1)$, so, $v=v_1$).
Fibres of $\vcenter{\xymatrix{F^*\bold{E} \ar[d]_-{p'} \\
\bold{A}}}$ are skeletal (assume, $(1_A,v):(A,E')@>\sim >>(A,E)$ is an iso in $(F^*\bold{E})_A$ then $v:E'@>\sim >>E$ is an iso in $\bold{E}_{F(A)}$, so, $E'=E$, and $(A,E')=(A,E)$). \hfill $\square $ 
\enddemo
\vskip 0.2cm
Partial cases of a pullback are 'fiber' and 'intersection of structures': 
\vskip 0.2cm
$\xymatrix{\bold{E}_B \ar[d] \ar[r] & \bold{E} \ar[d]^-{p} \\
1 \ar[r]_-{B} & \bold{B} }$
\hskip 2.5cm
$\xymatrix{\bold{E}_1\wedge \bold{E} \ar[d] \ar[r] \ar[dr]^-{\pi } & \bold{E} \ar[d]^-{p} \\
\bold{E}_1 \ar[r]_-{p_1} & \bold{B} }$

\vskip 0.2cm
{\bf Remark.} The notion of 'structure on objects' of a category was introduced in \cite{Kom} in order to deal with usual structures in Differential Geometry like smooth structures or fibre bundles on topological spaces. However, it turned out too weak (no inverse and direct images) and too strict (skeletal fibres) simoultaneously. Appropriate framework was created with theory of (co)fibrations. Nevertheless, a weaker notion of almost structure emphasizes direct connection with the main structure, which is the case of importance especially when almost structure is introduced on hom-sets.
\vskip 0.2cm
\proclaim{\bf Proposition 2.3} For each structure $\vcenter{\xymatrix{\bold{E} \ar[d]_-{p} \\ \bold{B}}}$ of type $\bold{E}$ on objects of category $\bold{B}$
\item{$\bullet $} there is an embedding
$i_p:\vcenter{\xymatrix{\bold{E} \ar[d]_-{p} \\ \bold{B}}} \hookrightarrow  \bold{ \ Set}^{\bold{E}^{op}}:\cases \dbinom E{p(E)}\mapsto p(\bold{E}(-,E)) & \text{ on objects}\\
\dbinom v{p(v)}\mapsto p(\bold{E}(-,v)) & \text{ on arrows}
\endcases $
\vskip 0.1cm
\item{$\bullet $} $p(\bold{E}(-,E))\hookrightarrow \bold{B}(p(-),p(E)):\bold{E}^{op}\to \bold{Set}$ (hom-subfunctor)
\endproclaim
\demo{Proof} Functoriality is obvious. Injectivity follows from Lemma 2.1.1.1.  \hfill $\square $
\enddemo
\vskip 0.2cm
{\bf Remark.} It means that every structure $\bold{E}$ on objects in $\bold{B}$ is {\bf faithfully} representable by a specific subcategory of $\bold{B}$-hom-subfunctors 
(in which sufficient to take arrows of only simple type $f\circ -\, $).
\vskip 0.2cm
The reasonable question is if can an object $E\in Ob\, \bold{E}$ be recovered from a functor $F\hookrightarrow \bold{B}(p(-),p(E))$? If can it is unique but the answer is no, in general.
Even when it is impossible subfunctor $F\hookrightarrow \bold{B}(p(-),p(E))$ behaves like an object in $\bold{E}$.
\vskip 0.2cm
{\bf Definition 2.2.}
\vskip 0.2cm
\item{$\bullet $} Arbitrary subfunctor $F\hookrightarrow \bold{B}(p(-),B):\bold{E}^{op}\to \bold{Set}$ is called {\bf almost-$\bold{E}$ structure} over object $B\in Ob\, \bold{B}$.
\item{$\bullet $} Category $\vcenter{\xymatrix{A\bold{E} \ar[d] \\ \bold{B}}}$ with objects $\dbinom FB, \ B\in Ob\, \bold{B}, \ F\hookrightarrow \bold{B}(p(-),B)$ and morphisms $\dbinom {f\circ -}f\equiv \dbinom {\bold{B}(p(-),f)}f, \ f:B\to B'$ is called a {\bf category of almost-$\bold{E}$ structures over $\bold{B}$}.
\vskip 0.1cm
\item{$\bullet $} {\bf Almost-$\bold{E}$ costructure over $B\in Ob\, \bold{B}$} is a subfunctor $F'\hookrightarrow \bold{B}(B,p(-)):\bold{E}\to \bold{Set}$ [almost costructures are not dual to almost structures, they behave all together in a covariant way].
\vskip 0.1cm
\item{$\bullet $} Category $\vcenter{\xymatrix{A'\bold{E} \ar[d] \\ \bold{B}}}$ with objects $\dbinom {F'}B, \ B\in Ob\, \bold{B}, \ F'\hookrightarrow \bold{B}(B,p(-))$ and morphisms $f:\dbinom {F'}B\to \dbinom {F'_1}{B_1}$,  if $f:B\to B_1\in Ar\, \bold{B}$ and $\forall E_1\in Ob\, \bold{E}, \ \forall g\in \bold{B}(B_1,p(E_1)). \ g\circ f\in \bold{B}(B,p(E_1))$,\linebreak
\vskip 0.1cm
is called a {\bf category of almost-$\bold{E}$ costructures over $\bold{B}$}. \hfill $\square $
\vskip 0.3cm
\centerline{\bf Example}
\vskip 0.3cm
Take $\bold{Poly}(E_1,E_2,...,E_n;-)\hookrightarrow \bold{Set}(p(E_1\times \cdots \times E_n),p(-)):\bold{Vect}\to \bold{Set}$, 
subfunctor of polylinear maps. Then $\bold{Poly}(+,+,...,+;-):\bold{Vect}^n\to \bold{Set}^{\bold{Vect}}$ determines a subcategory of 
almost-$\bold{Vect}$ costructures over $(p\circ (\times ^n))(\bold{Vect}^n)\hookrightarrow \bold{Set}$ ($p:\bold{Vect}\to \bold{Set}$ is forgetful).
\vskip 0.4cm

\proclaim{\bf Proposition 2.4} For a structure of type $\bold{E}$ on $\bold{B}$ \ $\vcenter{\xymatrix{\bold{E} \ar[d]_-{p} \\ \bold{B}}}$  
\newline
\vskip 0cm \hskip 0.5cm $\bullet $ $\vcenter{\xymatrix{A\bold{E} \ar[d] \\ \bold{B}}}$ is a fibration; \hskip 1cm
$\bullet $ $\vcenter{\xymatrix{\bold{E} \ar[d]_-{p} \\ \bold{B}}}\hookrightarrow \vcenter{\xymatrix{A\bold{E} \ar[d] \\ \bold{B}}}$
is a subcategory
\endproclaim
\demo{Proof} \item{$\bullet $} If \ $\dbinom FB\in Ob\, \vcenter{\xymatrix{A\bold{E} \ar[d] \\ \bold{B}}}$ and \  
$f:B'\to B$ \ take \ \newline
$f^*F:=\{g:p(X)\to B'\ |\ X\in Ob\, \bold{E},\ f\circ g\in F(X)\subset \bold{B}(p(X),B) \}$. Then $f^*F\hookrightarrow \bold{B}(p(-),B')$
is a subfunctor, and $\dbinom {f^*F}{B'}@>\tbinom {f\circ -}{f}>>\dbinom FB$ is cartesian over $f$ \hskip 0.5cm $\vcenter{\xymatrix{F'' \ar@{-->}[dr]_-{k\circ -} \ar[drr]^-{f\circ k\circ -} && \\
B'' \ar[drr]^-{f\circ k} \ar[dr]_-{k} & f^*F \ar[r]^-{f\circ -} & F \\
 & B' \ar[r]^-{f} & B
}}$
\item{$\bullet $} Assignment $\cases \dbinom E{p(E)}\mapsto p(\bold{E}(-,E))\hookrightarrow \bold{B}(p(-),p(E)) & \text{ on objects}\\
\dbinom v{p(v)}\mapsto p(\bold{E}(-,v))=\bold{B}(p(-),p(v)) & \text{ on arrows}
\endcases $
\vskip 0.3cm
gives the required embedding. \hfill $\square $   
\enddemo

\head {\bf 3. Enrichment with generalized elements in hom-sets} \endhead

{\bf Definition 3.1.} $1$-category $\bold{C}$ is {\bf enriched in tensor category}
$(\Cal V,I,\otimes )$ \cite{Kel, Bor2} if
\item{$\bullet $} $\forall x,y\in \text{\rm Ob}\, \bold{C}$ \ $\bold{C}(x,y)\in \text{\rm Ob}\Cal V$
\item{$\bullet $} $\forall x,y,z\in \text{\rm Ob}\, \bold{C}$ \ $\mu _{x,y,z}:\bold{C}(y,z)\otimes \bold{C}(x,y)\to \bold{C}(x,z)\in \text{\rm Ar}\Cal V$
\item{$\bullet $} $\forall x,y,z,w\in \text{\rm Ob}\, \bold{C}$ \newline
$$\xymatrix{(\bold{C}(z,w)\otimes \bold{C}(y,z))\otimes \bold{C}(x,y) \ar[r]^-{\sim } \ar[d]_-{\mu \otimes 1} & \bold{C}(z,w)\otimes (\bold{C}(y,z)\otimes \bold{C}(x,y)) \ar[r]^-{1\otimes \mu } & \bold{C}(z,w)\otimes \bold{C}(x,z) \ar[d]^-{\mu } \\
\bold{C}(w,y)\otimes \bold{C}(x,y) \ar[rr]_-{\mu } && \bold{C}(x,w)}$$
\item{$\bullet $} $\forall x\in \text{\rm Ob}\, \bold{C}$ $\exists \, u_x:I\to \bold{C}(x,x)\in \text{\rm Ar}\Cal V$ such that \newline
\vskip 0.1cm
\hbox{\hbox{$\xymatrix{\bold{C}(x,y)  & I\otimes \bold{C}(x,y) \ar[d]^-{u_x\otimes 1} \ar[l]_-{\sim } \\
 & \bold{C}(x,x)\otimes \bold{C}(x,y) \ar[ul]^-{\mu }
}$}
\hskip 2cm
\hbox{$\xymatrix{\bold{C}(z,x) & \bold{C}(z,x)\otimes I \ar[d]^-{1\otimes u_x} \ar[l]_-{\sim } \\
 & \bold{C}(z,x)\otimes \bold{C}(x,x) \ar[ul]^-{\mu }
}$}} \hfill $\square $

Generalized elements in hom-sets are just parametrized families of arrows in the same way as continuous or smooth families of maps. They form almost-$\bold{C}$ structure on hom-sets of category $\bold{C}$ \cite{Kom}.

\vskip 0.2cm
{\bf Definition 3.2.} Assume, $\bold{C}$ has binary products and $|\hskip 0.05cm\text{-}\hskip 0.05cm|:\bold{C}\to \bold{Set}$ is a faithful functor. Then $\bold{Set}$-map $f:|Z|\to \bold{C}(X,Y)$, $Z\in \text{\rm Ob}\, \bold{C}$ is called a {\bf generalized} (or acceptable) {\bf element} of $\bold{C}(X,Y)$ with domain $|Z|$ if arrow $\tilde f\circ \gamma _{Z}$ can be lifted to $\bold{C}$
$$\xymatrix{ \bold{C}(X,Y)\times |X| \ar[r]^-{ev} & |Y|\\
|Z|\times |X| \ar[u]^-{f\times 1} \ar[ur]_-{\tilde f} & \\
|Z\times X| \ar[u]^-{\gamma _{Z}} &
}$$
where: $\gamma _{Z}$ is the mediating arrow to product, $ev(g,x):=|g|(x)$ \hfill $\square $

\vskip 0.2cm
Denote $\Cal G(|Z|,\bold{C}(X,Y))\hookrightarrow \bold{Set}(|Z|,\bold{C}(X,Y))$ subset of
generalized elements of $\bold{C}(X,Y)$ with domain $|Z|$.

\proclaim{\bf Proposition 3.1} Assignment $Z\mapsto \Cal G(|Z|,\bold{C}(X,Y))$ is extendable to a functor \newline
$\Cal G(|\hskip 0.05cm\text{-}\hskip 0.05cm|,\bold{C}(X,Y)):\bold{C}^{op}\to \bold{Set}$.
\endproclaim
\demo{Proof} Assume, $\alpha :Z'\to Z$ is an arrow, $f\in \Cal G(|Z|,\bold{C}(X,Y))$ is a generalized element. We need to show that $f\circ \alpha \in \Cal G(|Z'|,\bold{C}(X,Y))$ is a generalized element as well, i.e., that $\exists h:Z'\times X\to Y$ such that $|h|=ev\circ (f\times 1)\circ (\alpha \times 1)\circ \gamma _{Z'}=\tilde f\circ (\alpha \times 1)\circ \gamma _{Z'}$.

Since $\gamma _{Z}:|Z\times X|\to |Z|\times |X|$ is natural in $Z$ \hskip 0.5cm
$\xymatrix{|Z'\times X| \ar[r]^-{\gamma _{Z'}} \ar[d]_-{|\alpha \times 1|} & |Z'|\times |X| \ar[d]^-{|\alpha |\times |1|} & \\
|Z\times X| \ar[r]_-{\gamma _Z} & |Z|\times |X| \ar[r]_-{\tilde f} & |Y|
}$\newline
the requirement will be $\exists h:Z'\times X\to Y$ such that $|h|=\tilde f\circ \gamma _Z\circ |\alpha \times 1|$. By assumption on $f$, $\exists g:Z\times X\to Y$ such that $|g|=\tilde f\circ \gamma _Z$. So, take $h:=g\circ (\alpha \times 1)$. \hfill $\square $
\enddemo

\proclaim{\bf Proposition 3.2} If $|\hskip 0.05cm\text{-}\hskip 0.05cm|:\bold{C}\to \bold{Set}$ is a faithful functor which preserves
binary products, then category $\bold{C}$ is enriched with generalized elements in presheaves category $\bold{Set}^{\bold{C}^{op}}$.
\endproclaim
\demo{Proof} \item{$\bullet $} $\forall X,Y\in Ob\, \bold{C}$ \ $\bigl (\Cal G(|\hskip 0.05cm\text{-}\hskip 0.05cm|,\bold{C}(X,Y)):\bold{C}^{op}\to \bold{Set}\bigr )\in Ob\bigl (\bold{Set}^{\bold{C}^{op}}\bigr )$
\item{$\bullet $} $\forall X,Y,Z\in Ob\, \bold{C}$ \ take \ $\mu _{X,Y,Z}:\Cal G(|\hskip 0.05cm\text{-}\hskip 0.05cm|,\bold{C}(Y,Z))\otimes \Cal G(|\hskip 0.05cm\text{-}\hskip 0.05cm|,\bold{C}(X,Y))\to \Cal G(|\hskip 0.05cm\text{-}\hskip 0.05cm|,\bold{C}(X,Z))$ such that
$\forall W\in Ob\, \bold{C} \ \ \mu _{X,Y,Z;W}(f,g):=\mu _{X,Y,Z}^{\bold{C}}\ \circ <f,g>$ where: $\mu ^{\bold{C}}_{X,Y,Z}:\bold{C}(Y,Z)\times \bold{C}(X,Y)\to \bold{C}(X,Z)$ is the composite in $\bold{C}$, \ $<f,g>:|W|\to \bold{C}(Y,Z)\times \bold{C}(X,Y)$ is the mediating arrow to product.
$\mu _{X,Y,Z;W}$ is natural in $W$ since $(\mu _{X,Y,Z}^{\bold{C}}\ \circ <f,g>)\circ |h|=\mu _{X,Y,Z}^{\bold{C}}\ \circ <f\circ |h|,g\circ |h|>$ for $h:W'\to W$.

Why $ev\circ ((\mu ^{\bold{C}}\ \circ <f,g>)\times 1)\circ \gamma $ can be lifted to $\bold{C}$? By condition,

\vskip 0.1cm
\hbox{
\hbox{$$\xymatrix{\bold{C}(Y,Z)\times |Y| \ar[r]^-{ev} & |Z| \\
|W|\times |Y| \ar[u]^-{f\times 1} \ar[ur]_-{\tilde f} & \\
|W\times Y| \ar[u]^-{\gamma } \ar[uur]_-{\text{lifted}}&
}$$}
\hskip 2cm
\hbox{$$\xymatrix{\bold{C}(X,Y)\times |X| \ar[r]^-{ev} & |Y| \\
|W|\times |X| \ar[u]^-{g\times 1} \ar[ur]_-{\tilde g} & \\
|W\times X| \ar[u]^-{\gamma } \ar[uur]_-{\text{lifted}}&
}$$}}
\vskip 0.2cm
Sufficient to take $\gamma =1$.
$$\xymatrix{\bold{C}(X,Z)\times |X| \ar@{..>}[rr]^-{ev} && |Z| \\
\bold{C}(Y,Z)\times \bold{C}(X,Y)\times |X| \ar[rr]^-{1_{\bold{C}(Y,Z)}\times ev} \ar@{..>}[u]^-{\mu ^{\bold{C}}\times 1_{|X|}} &&  \bold{C}(Y,Z)\times |Y| \ar[u]_-{ev} \\
 && \\
 && |W|\times |W|\times |X| \ar[uull]|-{f\times g\times 1_{|X|}} \ar[uu]|-{f\times \tilde g=\atop {(f\times 1_{|Y|})\circ (1_{|W|}\times \tilde g)}} \ar@/_12ex/[uuu]|-{\tilde f\circ (1_{|W|}\times \tilde g)} \\
 && |W|\times |X| \ar@{..>}[uuull]^-{<f,g>\times 1_{|X|}} \ar[u]_-{<1_{|W|},1_{|W|}>\times 1_{|X|}} \\
 && |W\times X| \ar@{..>}[u]_-{\gamma =1}
}$$
The required dotted way $ev\circ ((\mu ^{\bold{C}}\ \circ <f,g>)\times 1_{|X|})\circ \gamma $ can be lifted to $\bold{C}$ since the most right way $\tilde f\circ (1_{|W|}\times \tilde g)\circ (<1_{|W|},1_{|W|}>\times 1_{|X|})\circ 1_{|W\times X|}$ can be lifted (for take liftings for $\tilde f,\tilde g$ and identities for identities).
\item{$\bullet $} (associativity) $\forall f,g,h$ \ such that \ $f:|W|\to \bold{C}(Z,Z'), \ g:|W|\to \bold{C}(Y,Z), \ h:|W|\to \bold{C}(X,Y)$ \newline
$\mu ^{\bold{C}}_{X,Z,Z'}\, \circ <f,\mu ^{\bold{C}}_{X,Y,Z}\, \circ <g,h>>=\mu ^{\bold{C}}_{X,Y,Z'}\, \circ <\mu ^{\bold{C}}_{Y,Z,Z'}\, \circ <f,g>,h>$ because the equality holds at each point $w\in |W|$.
\item{$\bullet $} (identities) $\forall X\in Ob\, \bold{C}$ \ take
$u_{X;W}:\bold{1}\to \bold{Set}(|W|,\bold{C}(X,X)):*\mapsto (w\mapsto 1_X)$. It is natural in $W$ and $\forall f,g$ such that $f:|W|\to \bold{C}(X,Y)$, \, $g:|W|\to \bold{C}(Z,X)$
the equalities hold \linebreak
$\mu ^{\bold{C}}\, \circ <f,u_{X;W}>(w)=\mu ^{\bold{C}}(f(w),1_X)=f(w)$, \
$\mu ^{\bold{C}}\, \circ <u_{X;W},g>(w)=\mu ^{\bold{C}}(1_X,g(w))=g(w)$ \, $w\in |W|$.                   \hfill $\square $
\enddemo

\proclaim{\bf Corollary {\rm (refinement of Proposition 3.2)}} Under the above assumptions ($\bold{C}$ has binary products, $|\hskip 0.05cm\text{-}\hskip 0.05cm|:\bold{C}\to \bold{Set}$
is a faithful functor preserving binary products) hom-sets of $\bold{C}$ enriched with almost-$\bold{C}$ structure of generalized elements.
\endproclaim
\demo\nofrills{Proof \usualspace} is immediate (because presheaves of generalized elements are of specific form\linebreak
$\Cal G(|\hskip 0.05cm\text{-}\hskip 0.05cm|,\bold{C}(X,Y))\hookrightarrow \bold{Set}(|\hskip 0.05cm\text{-}\hskip 0.05cm|,\bold{C}(X,Y))$ and $\mu $ is actually postcomposite $\mu ^{\bold{C}}\circ -$).    \hfill $\square $
\enddemo
\vskip 0.1cm
We mean further that $\bold{C}$ is $A\bold{C}$-category if this {\bf specific enrichment} with generalized elements is given. Moreover, we call $\bold{D}$ is $A\bold{C}$-category if it is enriched with presheaves of generalized elements with domains in $\bold{C}$.
\vskip 0.1cm

{\bf Remark.} All usual concrete categories, like $\bold{Top}, \ \bold{Grp}, \ \bold{Rng}$, etc. carry corresponding almost structures (which in some cases can be strict).
\vskip 0.35cm
\centerline{\bf Example}
\vskip 0.3cm

\proclaim{\bf Proposition 3.3} If $X$ is a locally compact topological space
(so that, family $\Cal T$ of topologies on $\bold{Top}(X,Y)$ for which evaluation map $ev:\bold{Top}(X,Y)\times |X|\to |Y|$ is continuous
 is not empty and contains minimal element, compact-open topology on $\bold{Top}(X,Y)$) then $\tau \in \Cal T$ is compact-open
 iff $\forall Z\in Ob\, \bold{Top}$ each generalized element $f:|Z|\to \bold{Top}(X,Y)$ is continuous.
\endproclaim
\demo{Proof} $"\Rightarrow "$ \ Regard the diagram
$$\xymatrix{\bold{Top}(X,Y)\times |X| \ar[rr]^-{ev} && |Y| \\
|Z|\times |X| \ar[u]^-{f\times 1} \ar[urr]_-{\tilde f} &&
}$$
We want to show that $f:|Z|\to \bold{Top}(X,Y)$ is continuous (with compact-open topology in $\bold{Top}(X,Y)$) if 
$\tilde f:|Z|\times |X|\to |Y|$ is continuous, i.e., that $\forall z\in |Z|$ \ $\forall $ (subbase) compact-open set $U^K\ni f(z)$ 
\ $\exists $ nbhd $V\ni z$ such that $f(V)\subset U^K$. It is equivalent that $\forall z\in |Z| \ \forall U^K\ni f(z) \ \exists V\ni z$ 
such that $\tilde f(V\times K)\subset U$. Since $\tilde f$ is continuous $\forall (z,x)\in \{z\}\times K$ and $\forall $ open 
$U\ni \tilde f(z,x)$ \ $\exists $ open nbhd $V_z\times W_x\ni (z,x)$ such that $\tilde f(V_z\times W_x)\subset U$. 
$\bigcup\limits _{x\in K}W_x\supset K$ (open cover). So, by compactness of $K$, $\exists W_{x_1},...,W_{x_n}$, such that 
$\bigcup\limits _{i=1}^nW_{x_i}\supset K$. Take $V:=\bigcap\limits _{i=1}^n(V_i)_z$, where $(V_i)_z$ corresponds to $W_{x_i}$ 
(i.e., $(V_i)_z$ is open, $(V_i)_z\ni z$, $\tilde f((V_i)_z\times W_{x_i})\subset U$). Then $\tilde f(V\times K)\subset U$. 
\vskip 0.2cm
$"\Leftarrow "$ Take $Z=\bold{Top}(X,Y)$ with compact-open topology. Take $\bold{Top}(X,Y)$ itself (on the top of the diagram) with 
non minimal $\tau \in \Cal T$, \ $f:=1\in Ar\, \bold{Set}$. Then $1:\bold{Top}(X,Y)\to \bold{Top}(X,Y)$
is an admissible generalized element, since $ev$ is continuous, but $1$ is not continuous.  
\hfill $\square $
\enddemo
\vskip 0.2cm
{\bf Remark.} Therefore, for locally compact space $X$ almost-$\bold{Top}$ structure $\Cal G(|Z|,\bold{Top}(X,Y))$ coincides with 
compact-open topology, i.e., is actually $\bold{Top}$ structure.

\vskip 0.3cm
If we agree that functor $\Cal G(|\hskip 0.05cm\text{-}\hskip 0.05cm|,\bold{C}(X,Y)):\bold{C}^{op}\to \bold{Set}$ reflects essential 
properties of $\bold{C}$-hom-sets we immediately get a unique (up to isomorphism) extension of each functor $F:\bold{C}\to \bold{D}$, i.e., 
deal with $\bold{C}$-hom-sets as with $\bold{C}$-objects. In this way, for example, tangent or jet functor can be introduced directly on 
$\bold{Aut}(X)$, $X\in Ob\, \bold{Diff}$ to give rise a calculus on $\bold{Aut}(X)$. Possibility of such an extension follows from the
fact that each presheaf is a certain (canonical) colimit of representables \cite{Mac, M-M}.

\proclaim{\bf Proposition 3.4} 
\item{$\bullet $} Yoneda embedding $y:\bold{C}\to \bold{Set}^{\bold{C}^{op}}$ is a universal cocompletion of 
$\bold{C}$, i.e., $\forall \, F:\bold{C}\to \bold{E}$, where $\bold{E}$ is cocomplete, $\exists \, ! \text{ (up to iso) cocontinuous } 
\bar F:\bold{Set}^{\bold{C}^{op}}\to \bold{E}$ such that \hskip 0.35cm
$\vcenter{\xymatrix{\bold{Set}^{\bold{C}^{op}} \ar@{-->}[r]^-{\bar F} & \bold{E} \\
\bold{C} \ar[u]^-{y} \ar[ur]_-{F} & }}$
\item{$\bullet $} $\bar F(P)=Colim(\int _{\bold{C}}P@>\pi >>\bold{C}@>F>>\bold{E})$, where $P\in Ob\, \bold{Set}^{\bold{C}^{op}}$, 
$\int _{\bold{C}}P$ is a category of elements of $P$, $\pi $ is the natural projection.
\item{$\bullet $} $\xymatrix{\bold{Cat} \ar@/^1.5ex/[rr]^-{\bold{Set}^{(-)^{op}}} & & \bold{Cocomp} \ar@/^1.5ex/[ll]^-{\text{forgetful}}_-{\bot }}$, 
adjunction between $\bold{Cat}$ and full subcategory of cocomplete categories $\bold{Cocomp}$ with Yoneda embedding $y_{\bold{C}}:\bold{C}\to \bold{Set}^{\bold{C}^{op}}$ as a unit.
\item{$\bullet $} Each functor $F:\bold{C}\to \bold{D}$ admits a unique (up to iso) cocontinuous extension $F:\bold{Set}^{\bold{C}^{op}}\to \bold{Set}^{\bold{D}^{op}}$
such that \hskip 0.6cm
$\xymatrix{\bold{Set}^{\bold{C}^{op}} \ar@{-->}[r]^-{F} & \bold{Set}^{\bold{D}^{op}}\\
\bold{C} \ar[u]^-{y_{\bold{C}}} \ar[r]^-{F} & \bold{D} \ar[u]_-{y_{\bold{D}}}}$
\endproclaim
\demo{Proof} \ See \cite{M-M}.  \hfill $\square $
\enddemo

For example, if $T:\bold{Diff}\to \bold{Diff}$ is a tangent functor, $\bold{Diff}(X,Y)$ is a presheaf on $\bold{Diff}$ 
(hom-set enriched as above) then $T(\bold{Diff}(X,Y))=\int _{\bold{Diff}}\bold{Diff}(X,Y)@>\pi>>\bold{Diff}@>T>>\bold{Diff}@>y>>\bold{Set}^{\bold{Diff}^{op}}$.

\vskip 0.1cm
\subhead 3.1. Tangent functor for smooth algebras\endsubhead

It is an example of dual (and invariant) construction for the main functor of Differential Geometry (which gives suggestion how it can be extended over spectra of commutative algebras).

Let $T:\bold{Diff}\to \bold{Diff}$ be {\bf tangent} functor on the category of real $\infty $-smooth manifolds. 
In local 
\vskip 0.2cm
coordinates it looks like $\cases X\to TX:(x^i)\to (x^i,\xi ^j) & X\in Ob\, \bold{Diff}\\
f\to Tf:(f^i(x))\to (f^i(x),\frac {\partial f^j}{\partial x^k}\, \xi ^k) & f\in Ar\, \bold{Diff}
\endcases $
\vskip 0.2cm
$\bold{Diff}\hookrightarrow \Bbb R\text{-}\bold{Alg}^{op}$ is a subcategory of the opposite of real commutative algebras. 
Working in $\bold{Diff}$ it is hard (if possible at all) to give coordinate-free characterization of $T$. The question is 
how it looks like in $\Bbb R\text{-}\bold{Alg}$?
\vskip 0.2cm
{\bf Definition 3.1.1.} Let $\Cal A\in Ob\, \Bbb R\text{-}\bold{Alg}$.  
\item{$\bullet $} $\rho :\Cal A\to \bold{Top}(\bold{Spec}_{\Bbb R}(\Cal A),\Bbb R)$
is called {\bf functional representation} homomorphism of $\Cal A$, where $\bold{Spec}_{\Bbb R}(\Cal A)=\Bbb R\text{-}\bold{Alg}(\Cal A,\Bbb R)$ with initial 
topology w.r.t. all functions $\rho (a), \ a \in \Cal A$, \ $\rho (a)(f):=ev(f,a):=|f|(a)$. 
\item{$\bullet $} $\Cal A$ is called {\bf smooth} if $\forall a_1,a_2,...,a_n\in \Cal A$ and 
$\forall f:\Bbb R^n\to \Bbb R\in C^{\infty }(\Bbb R^n)$ the composite \newline
$f\, \circ <\rho (a_1),\rho (a_2),...,\rho (a_n)>\ \in \ \text{Im}\, (\rho )$. \hfill $\square $
\vskip 0.2cm
Denote by $\Bbb R\text{-{\tt Sm}-}\bold{Alg}\hookrightarrow \Bbb R\text{-}\bold{Alg}$ full subcategory of smooth algebras.
\vskip 0.2cm
\proclaim{\bf Lemma 3.1.1} $\Bbb R\text{-{\tt Sm}-}\bold{Alg}\hookrightarrow \Bbb R\text{-}\bold{Alg}$ is a reflective subcategory, i.e., 
the inclusion has a left adjoint $\text{\tt Sm}:\Bbb R\text{-}\bold{Alg}\to \Bbb R\text{-{\tt Sm}-}\bold{Alg}$, {\bf smooth completion} of 
$\Bbb R$-algebras.
\endproclaim
\demo{Proof} Just take for each $\Bbb R$-algebra $\Cal A$ $\Bbb R$-algebra $\text{\tt Sm}(\Cal A)$ of all terms $\{f(a_1,...,a_n)\, |\, f:\Bbb R^n\to \Bbb R,\ a_1,...,a_n\in \Cal A\}$
(all smooth operations are admitted). Each morphism $f$ from an $\Bbb R$-algebra $\Cal A$ to a smooth alebra $\Cal B$ is uniquely 
extendable to $\tilde  f:\text{\tt Sm}(\Cal A)\to \Cal B$. \hfill $\square $
\enddemo

Let $\bold{Sym}\text{-}\bold{Alg}$ be a category of symmetric partial differential algebras. $Ob\, (\bold{Sym}\text{-}\bold{Alg})$ are 
graded commutative algebras over commutative $\Bbb R$-algebras with a differential $d:\Cal A^0\to \Cal A^1$ of degree $1$  
determined only on elements of degree $0$ ($d$ is $\Bbb R$-linear and satisfies Leibniz rule). $Ar\, (\bold{Sym}\text{-}\bold{Alg})$ are 
graded degree $0$ algebra homomorphisms which respect $d$. 

\proclaim{\bf Lemma 3.1.2} There is an adjunction $\vcenter{\xymatrix{\Bbb R\text{-}\bold{Alg} \ar@/^/[r]^-{\bold{Sym}}_-{\bot } & \bold{Sym}\text{-}\bold{Alg} \ar@/^/[l]^-{p_0}}}$, \newline 
where: $p_0$ is the projection onto $0$-component \ $\cases p_0(\Cal A):=\Cal A^0 & \\
p_0(\Cal A@>f>>\Cal B):=(\Cal A^0@>f^0>>\Cal B^0) & \endcases $, \vskip 0.2cm
$\bold{Sym}$ is taking symmetric algebra over module of differentials of the given algebra \vskip 0.2cm
$\cases \bold{Sym}(\Cal C):=\bold{Sym}(\Lambda ^1(\Cal C)) & \\
\bold{Sym}(\Cal C@>h>>\Cal D):=(\bold{Sym}(\Cal C)@>\tilde h>>\bold{Sym}(\Cal D)) & \\
\tilde h\bigl (\sum c_{i^1\dots i^k}(dc_1)^{i_1}\cdots (dc_k)^{i_k}\bigr ):=\sum h(c_{i^1\dots i^k})(dh(c_1))^{i_1}\cdots (dh(c_k))^{i_k} & \endcases $ \vskip -0.35cm \hfill $\square $
\endproclaim

\vskip 0.2cm
\proclaim{\bf Lemma 3.1.3} \ $\vcenter{\xymatrix{\Bbb R\text{-}\bold{Alg} \ar[r]^-{\text{\tt Sm}} \ar[dr]_-{\bold{Spec}_{\Bbb R}} & \Bbb R\text{-{\tt Sm}-}\bold{Alg} \ar[d]^-{\bold{Spec}_{\Bbb R}} \\
 & \bold{Top}}}$ \ (smooth completion does not change spectrum).
\endproclaim
\demo{Proof} $\forall \alpha :\Cal A\to \Bbb R$ $\exists ! \text{ extension } \tilde \alpha :\text{\tt Sm}(\Cal A)\to \Bbb R:f(a_1,\dots ,a_n)\mapsto f(\alpha (a_1),\dots ,\alpha (a_n))$. And conversely,
each such $\tilde \alpha $ is uniqely restricted to $\alpha $. Initial topology on $\Bbb R\text{-}\bold{Alg}((\text{\tt Sm})(\Cal A),\Bbb R)$ 
does not change because new functions are functionally dependent on old ones.    \hfill $\square $ 
\enddemo

{\bf Remark.} With Zarisski topology in spectra smooth completion yields the same set with a weaker topology. 
For \, $C^{\infty }(X),\, X\in Ob\, \bold{Diff}$ \, Zarisski and initial topologies coincide.  

\vskip 0.2cm
\proclaim{\bf Proposition 3.1.1} \,  
$\bullet $ Tangent functor $T:\Bbb R\text{-{\tt Sm}-}\bold{Alg}\to \Bbb R\text{-{\tt Sm}-}\bold{Alg}$ is equal to the composite
$\Bbb R\text{-{\tt Sm}-}\bold{Alg}\hookrightarrow \Bbb R\text{-}\bold{Alg}@>\bold{Sym}>> \bold{Sym}\text{-}\bold{Alg}@>U>> \Bbb R\text{-}\bold{Alg}@>\text{\tt Sm}>> \Bbb R\text{-{\tt Sm}-}\bold{Alg}$, 
where $U$ forgets differential $d$ and grading.
\item{$\bullet $} To canonical projection $\vcenter{\xymatrix{TX \ar[d]_-{p_X} \\ X}}$ there corresponds canonical embedding 
$\vcenter{\xymatrix{T(C^{\infty }(X))  \\ C^{\infty }(X) \ar[u]^-{i_{C^{\infty }(X)}}}}$.
\endproclaim
\demo{Proof} 
\item{$\bullet $} If $X\in Ob\, \bold{Diff} $ \ $TX\sim \bold{Spec}_{\Bbb R}(U\circ \bold{Sym}(C^{\infty }(X)))\sim \bold{Spec}_{\Bbb R}(\text{\tt Sm}\hskip 0.05cm\circ U\circ \bold{Sym}(C^{\infty }(X)))$.                                   
\item{$\bullet $} immediate.                 \hfill $\square $
\enddemo
\vskip 0.2cm
{\bf Remark.} It is reasonable to define $T$ on $\Bbb R\text{-}\bold{Alg}$ as $T:=U\hskip -0.05cm\circ \bold{Sym}$ and transfer 
it to spectra via duality $\xymatrix{\Bbb R\text{-}\bold{Alg}^{op} \ar@/^/[r]_-{\bot }^-{F} & \bold{Spec}_{\Bbb R} \ar@/^/[l]^-{G}}$ (as $F\circ T^{op}\circ G$).

\head {\bf 4. General manifolds} \endhead

{\bf Definition 4.1.} Functor $\vcenter{\xymatrix{\bold{E} \ar[d]_-{F} \\ \bold{B} }}$ is called a {\bf fibration with respect to} 
class of arrows $\Cal C\subset Ar\, \bold{B}$ if $\forall f:B'\to F(E)\in \Cal C$ \ $\exists \tilde f:E'\to E\in Ar\, \bold{E}$ such that 
$\tilde f$ is over $f$  and $\tilde f$ is cartesian.   \hfill $\square $
\vskip 0.2cm
{\bf Definition 4.2.} \cite{M-M} 
\vskip 0.1cm
\item{$\bullet $} Grothendieck {\bf pretopology} $\tau _0$ on a category $\bold{B}$ with pullbacks is a family of {\bf coverings}
$\tau _{0B}$ for each object $B\in Ob\, \bold{B}$ (elements of a covering are just arrows with codomain $B$) such that
\itemitem{$\bullet $} if $f:B'@>\sim >>B$ is an iso then $\{f\}\in \tau _{0B}$ is an one-element covering
\itemitem{$\bullet $} if $g:B''\to B$ is an arrow and $\goth c\in \tau _{0B}$ then pullback family $\text{\bf plbk}_g(\goth c)\in \tau _{0B''}$
\itemitem{$\bullet $} (coverings are composable) if $\goth c\in \tau _{0B}$ and $\forall B'\in d(\goth c)$ there is given a covering $\goth c_{B'}\in \tau _{0B'}$
then \ $\goth c\, \circ \hskip -0.1cm\bigcup\limits _{B'\in d(\goth c)}\hskip -0.1cm\goth c_{B'}\in \tau _{0B}$
\item{$\bullet $} Grothendieck {\bf topology} $\tau $ on a category $\bold{B}$ (not necessarily with pullbacks) is a family of hom-subfunctors $\tau _B$ 
for each object $B\in Ob\, \bold{B}$ such that
\itemitem{$\bullet $} $\bold{B}(-,B)\in \tau_B$
\itemitem{$\bullet $} if $f:B'\to B$ and $\goth t\in \tau _B$ then the {\bf inverse image} $\bigl (f^*(\goth t):X\mapsto (f^*\goth t)(X,B')\subset \bold{B}(X,B')\bigr )\in \tau _{B'}$ \ ($h\in f^*(\goth t)(X,B')$ 
iff $f\circ h\in \goth t(X,B)$)
\itemitem{$\bullet $} if $\goth s\hookrightarrow \bold{B}(-,B)$ is any hom-subfunctor such that $\forall f:B'\to B\ f^*(\goth s)\in \tau _{B'}$ 
then $\goth s\in \tau _B$    \hfill $\square $
\vskip 0.2cm
Every topology is a pretopology if $\bold{B}$ has pullbacks, and every pretopology generates a topology \cite{M-M}. Category $\bold{B}$
with (pre)topology is called {\bf site} $(\bold{B},\tau _{(0)})$.
\vskip 0.2cm
{\bf Definition 4.3.} \cite{Kom} Functor $\vcenter{\xymatrix{\bold{E} \ar[d]_-{F} \\ (\bold{B},\tau _0)}}$ is called {\bf local} if it is a fibration
with respect to all elements of all coverings \ $\bigcup \hskip -0.1cm\bigcup\limits _{B\in Ob\, \bold{B}}\hskip -0.1cm\tau _{0B}$. \hfill $\square $

{\bf Definition 4.4.} For a given functor to a site \hskip -0.0cm$\vcenter{\xymatrix{\bold{E} \ar[d]_-{F} \\ (\bold{B},\tau _0)}}$ \hskip -0.05cmsmallest local functor \hskip -0.05cm 
$\vcenter{\xymatrix{\bold{E}\text{-}\bold{Man} \ar[d]_-{p_F} & \bold{E} \ar@{_{(}->}[l] \ar[dl]^-{F}\\
(\bold{B},\tau _0) & }}$ is called {\bf $\bold{E}$-manifold structure over $\bold{B}$}. It means
\vskip 0.1cm
\item{$\bullet $} $\forall X\in Ob\, \bold{E}\text{-}\bold{Man}$ $\exists $ covering $\goth c_{p_F(X)}=\{i:B_i\to p_F(X)\}_{i\in I}\in \tau _{0p_F(X)}$ 
such that there are inverse images $i^*(X)\in Ob\, \bold{E}$ (i.e., $\bold{E}$ contains isomorphic representatives of the inverse images)
\vskip 0.1cm
\item{$\bullet $} $\forall f:X'\to X\in Ar\, \bold{E}\text{-}\bold{Man}$ \, $\exists $ \, coverings $\goth c'_{p_F(X')}=\{i':B'_{i'}\to p_F(X')\}\in \tau _{0p_F(X')}$ 
and $\goth c_{p_F(X)}=\{i:B_{i}\to p_F(X)\}\in \tau _{0p_F(X)}$ such that $\forall i\in \goth c_{p_F(X)}$ $\exists i'\in \goth c'_{p_F(X')}$ 
such that $\vcenter{\xymatrix{B'_{i'} \ar[r]^-{\exists \, \varphi } \ar[d]_-{i'} & B_i \ar[d]^-{i}\\
p_F(X') \ar[r]_-{p_F(f)} & p_F(X)}}$ and over it $\vcenter{\xymatrix{i^{'*}(X') \ar@{-->}[r]^-{\exists !\, \Phi } \ar[d]_-{\tilde {i'}} & i^{*}(X) \ar[d]^-{\tilde i}\\
X' \ar[r]_-{f} & X}}$ \newline
$\tilde {i'}, \tilde i$ are cartesian, \ $p_F(\Phi )=\varphi $, \ $\Phi \in Ar\, \bold{E}$ (arrows are locally in $\bold{E}$)
\vskip 0.1cm
\item{$\bullet $} $\bold{E}\text{-}\bold{Man}$ is maximal with respect to two above properties \hfill $\square $
\vskip 0.2cm
\centerline{\bf Examples}
\vskip 0.2cm
\item{1.} $\bold{Set}$ as a manifold structure $\vcenter{\xymatrix{\bold{Set} \ar[d]_-{1} & \bold{\, Set}_{inj} \ar@{_{(}->}[l] \ar@{_{(}->}[dl]\\
(\bold{Set},\tau _0) & }}$ \ where $\bold{Set}_{inj}$ is a category of sets with injective maps only, $\tau _0$ is a pretopology 
consisting of all families of injective maps with common codomains.
\item{2.} {\bf Differentiable} manifolds \hskip 0.5cm $\vcenter{\xymatrix{C^r\text{-}\bold{Man}_k \ar@{_{(}->}[d] & \text{\tt \, Triv}C^r\text{-}\bold{Man}_k \ar@{_{(}->}[l] \ar@{_{(}->}@<0.5ex>[dl] \\
(\bold{Top},\tau _0) & }}$ \hskip 0.5cm where $k=\Bbb R$ or $\Bbb C$, \newline
$\tau _0$ consists of all open coverings, \ $r=0,1,..., \infty $ (or $\omega $
for complex manifolds), \vskip 0.2cm
$\cases Ob\, (\text{\tt Triv}C^r\text{-}\bold{Man}_k)=\{k^0,k^1,..., k^n,...\} & \\ Ar\, (\text{\tt Triv}C^r\text{-}\bold{Man}_k)=C^r\text{-maps} & \endcases $
\vskip 0.2cm
\item{3.} Locally trivial {\bf fibre bundles} \hskip 0.25cm$\vcenter{\xymatrix{\bold{Bn}\, (\Cal E,p) \ar[d] & \bold{Bn}_0\, (\Cal E,p) \ar@{_{(}->}[l] \ar[dl] \\
(\bold{Man}^{\to },\tau _0) & }}$ \hskip 0.25cm(see 5)
\vskip 0.2cm
\item{4.} {\bf Foliations} over $\bold{Man}$ \hskip 0.25cm$\vcenter{\xymatrix{\bold{Fol}\, (\Cal E,p) \ar[d] & \bold{Fol}_0\, (\Cal E,p) \ar@{_{(}->}[l] \ar[dl] \\
(\bold{Man},\tau _0) & }}$ \hskip 0.25cm associated to $A\bold{Man}$-functor sequence $\Cal E@>p>>\bold{Man}@>\pi >>\bold{Top}$ \ 
(see 5), \ where:
$\bold{Man}$ is a category of manifolds (of type $\Cal E'$) over $\bold{Top}$, $\tau _0$ are all 'open coverings' of objects in $\bold{Man}$, 
$\bold{Fol}_0\, (\Cal E,p)\equiv \bold{Bn}_0\, (\Cal E,p)$ is a category of trivial foliations ('direct products') with leaves in $\Cal E$, 
projection functor $\vcenter{\xymatrix{\bold{Fol}_0(\Cal E,p) \ar[d] \\
(\bold{Man},\tau _0)}}$ is the first (top) component of projection $\vcenter{\xymatrix{\bold{Bn}_0\, (\Cal E,p) \ar[d] \\
(\bold{Man}^{\to },\tau _0)}}$ ($\tau _0$ are different in these two cases and corresponding categories of manifolds are glued differently).
\item{5.} $\bold{E}$-manifolds over $\bold{Top}$. Let $\vcenter{\xymatrix{\bold{E} \ar[d]_-{p} \\ (\bold{Top},\tau _0)}}$ be a local structure on $\bold{Top}$ with $\tau _0$, all open coverings.
\itemitem{$\bullet $} {\bf Local $\bold{E}$-map} on a topological space $X$ is a pair \ $\dbinom EU\in Ob\vcenter{\xymatrix{\bold{E} \ar[d] \\ (\bold{Top},\tau _0)}}$, \hskip 0.5cm$U$ is open.
\itemitem{$\bullet $} Family $\Bigl \{ \dbinom {E_i}{U_i}\Bigr \} _{i\in I}$ is {\bf compatible} iff $\forall (i,j)\in I^2$ \, $E_i|_{U_i\cap U_j}@>\sim >\varphi >U_j|_{U_i\cap U_j}$, $\varphi $ is a vertical iso.
\itemitem{$\bullet $} {\bf $\bold{E}$-atlas} $\Cal A$ on $X$ is a compatible family $\Bigl \{ \dbinom {E_i}{U_i}\Bigr \} _{i\in I}$ such that $\bigcup\limits _{i\in I}U_i=X$.
\itemitem{$\bullet $} Two $\bold{E}$-atlases $\Cal A$ and $\Cal A'$ are {\bf equivalent} iff $\Cal A\cup \Cal A'$ is still an $\bold{E}$-atlas on $X$ (so, there exist maximal atlases, call them $\Cal A_{max}$, $\Cal B_{max}$, etc.).
\itemitem{$\bullet $} The above 'equivalence' on atlases is not transitive in general. 
So, there can be different maximal atlases containing a given one. But, it is transitive if $\forall \, \dbinom EU, \dbinom {E'}{U}\in Ob\, \vcenter{\xymatrix{\bold{E} \ar[d]_-{p} \\ \bold{Top}}}$ 
and $\forall \text{\rm \ open covering } \bigcup\limits _{i\in I}U_i\supset U$ \ $E|_{U_i}@>\sim >\text{\rm vert}>E'|_{U_i}$ (for all $i\in I$)
implies $E@>\sim >\text{\rm vert}>E'$. 
\itemitem{$\bullet $} Topological space $X$ together with an $\bold{E}$-atlas $\Cal A$ on it is called {\bf $\bold{E}$-manifold}, i.e., $(X,\Cal A)\in Ob\, \bold{E}\text{-}\bold{Man}$.
\itemitem{$\bullet $} An {\bf arrow} in $\bold{E}\text{-}\bold{Man}$ is $f:(X,\Cal A)\to (Y,\Cal B)$ such that $f:X\to Y$ is 
continuous and $\forall \, \dbinom EU\in \Cal A, \dbinom {E'}{V}\in \Cal B$ if $U\cap f^{-1}(V)\ne \emptyset $ then $f|_{U\cap f^{-1}(V)}:U\cap f^{-1}(V)\to V$ 
admits (unique) lifting $\bar f|_{U\cap f^{-1}(V)}:E|_{U\cap f^{-1}(V)}\to E'\in Ar\, \bold{E}$.

\head {\bf 5. Fibre bundles} \endhead

Locally trivial fibre bundles give an important example of general manifolds over $\bold{Man}^{\to }$ \cite{Kom}.
\vskip 0.2cm
{\bf Definition 5.1.} Category of {\bf trivial fibre bundles} $\bold{Bn}_0\, (\Cal E,p)$ over $\bold{Man}^{\to }$ with {\bf typical fibres} in a category $\Cal E$ consists of the following data
\item{$\bullet $} $\Cal {E}@>p>>\bold{Man}@>\pi >>\bold{Top}$, \ where \ $\Cal E$ and $\bold{Man}$ are $A\bold{Man}$-categories, $p$ is $A\bold{Man}$-functor, $\pi $ preserves binary products [i.e., $\bold{Man}$ is enriched in $\bold{Set}^{\bold{Man}^{op}}$ with presheaves of generalized elements $\Cal G(|-|,\bold{Man}(A,A'))$ for each hom-set $\bold{Man}(A,A')$, $\Cal E$ is enriched in $\bold{Set}^{\bold{Man}^{op}}$ with subfunctors $\Cal H(|-|,\Cal E(E,E'))\hookrightarrow \bold{Set}(|-|,\Cal E(E,E'))$ for each hom-set $\Cal E(E,E')$, $p_{E,E';X}:\Cal H(|X|,\Cal E(E.E'))\to \Cal G(|X|,\bold{Man}(p(E),p(E'))):f\mapsto p_{E,E'}\circ f$ is natural in $X\in Ob\, \bold{Man}$, $p_{E,E'}:\Cal E(E,E')\to \bold{Man}(p(E),p(E'))$ is the restriction of functor $p$ on the hom-set]
\item{$\bullet $} $Ob\, \bold{Bn}_0\, (\Cal E,p):=\{(X,E)\, |\, X\in Ob\, \bold{Man},E\in Ob\, \Cal E\}$;\newline 
$Ar\, \bold{Bn}_0\, (\Cal E,p):=\{(X,E)@>(f,\Phi )>> (X',E')\, |\, f:X\to X',\ \Phi \in  \Cal H(|X|,\Cal E(E,E'))\}$ 
\vskip 0.1cm
\item{$\bullet $} functor $\vcenter{\xymatrix{\bold{Bn}_0\, (\Cal E,p) \ar[d]_-{p_0} \\ \bold{Man}^{\to }}}:\cases (X,E)\mapsto \hskip 1cmX\times p(E)@>p_1>>X & \hskip -0.1cm(X,E)\in Ob\, \bold{Bn}_0\, (\Cal E,p) \hskip -0.1cm\\
(f,\Phi )\mapsto \vcenter{\xymatrix{X\times p(E) \ar[rr]^-{<f\circ p_1,\phi >} \ar[d]_-{p_1} && X'\times p(E') \ar[d]^-{p_1}\\
X \ar[rr]_-{f} && X'}} & \hskip -0.05cm(f,\Phi )\in Ar\, \bold{Bn}_0\, (\Cal E,p)\hskip -0.1cm
\endcases $
\vskip 0.1cm
where $\phi :=ev\circ ((p_{E,E'}\circ \Phi )\times 1_{|p(E)|})$, \ $p_{E,E'}\circ \Phi \in \Cal G(|X|,\bold{Man}(p(E),p(E')))$
\hfill $\square $
\vskip 0.2cm
{\bf Definition 5.2.} Category of {\bf locally trivial fibre bundles} $\bold{Bn}\, (\Cal E,p)$ over site $(\bold{Man^{\to }},\tau _0)$, where $\tau _0$ is pullbacks of all open coverings of codomains 
(i.e., if $q:Y\to X\in Ob\, \bold{Man}^{\to }$ and $\{U_i\}_{i\in I}$ is an open covering of $X$ then $\{q^{-1}(U_i)@>q|_{U_i}>> U_i\}_{i\in I}$ is a covering of $q$), 
is a manifold structure of type $\bold{Bn}_0\, (\Cal E,p)$ over $\bold{Man}^{\to }$.
\hfill $\square $
\vskip 0.2cm
A usual way of construction of new fibre bundles from old ones is by fibrewise operations.
Let $\Cal E@>p>>\bold{Man}@>\pi >>\bold{Top}$ \, and \, $\Cal E'@>p'>>\bold{Man}@>\pi >>\bold{Top}$ \, be two sequences generating 
categories of fibre bundles \, $\bold{Bn}\, (\Cal E,p)$ \, and \, $\bold{Bn}\, (\Cal E',p')$ \, of types \, $\Cal E$ \, and \, $\Cal E'$, respectively, 
\, $F:\Cal E\to \Cal E'$ \, be an $A\bold{Man}$-functor. Then there exists a corresponding functor \, $\bold{Bn}\, (F):\bold{Bn}\, (\Cal E,p)\to \bold{Bn}\, (\Cal E',p')$.

\vskip 0.1cm
Denote by $A\bold{Man}\text{-}\bold{CAT}$ an $1$-category such that 
\vskip 0.1cm
$\cases \hskip -0.1cmOb\, (A\bold{Man}\text{-}\bold{CAT})\ni (\Cal E,p), \ \text{ if } \Cal E \text{ is } A\bold{Man}\text{-category}, \ p:\Cal E\to \bold{Man} \text{ is } A\bold{Man}\text{-functor} & \\ 
\hskip -0.1cmAr\, (A\bold{Man}\text{-}\bold{CAT}) \ni (F:(\Cal E,p)\to (\Cal E',p')), \ \text{ if } F:\Cal E\to \Cal E' \text{ is } A\bold{Man}\text{-functors} & \endcases $
\vskip 0.1cm
and by $\bold{Bn}_0$ and $\bold{Bn}$ subcategories of $1\text{-}\bold{CAT}$ consisting of categories of 
trivial and locally trivial fibre bundles with fibres of a fixed type (i.e., of categories like $\bold{Bn}_0\, (\Cal E,p)$ and $\bold{Bn}\, (\Cal E,p)$) and functors preserving atlases as arrows (see 2.5, remarks). 
Of course, $\bold{Bn}_0\, (\Cal E,p)\hookrightarrow \bold{Bn}\, (\Cal E,p)$.

\proclaim {\bf Proposition 5.1} There are functors \ $\bold{Bn}_0\, (-):A\bold{Man}\text{-}\bold{CAT}\to \bold{Bn}_0\hookrightarrow 1\text{-}\bold{CAT}:$ 
\vskip 0.15cm
$\cases  \hskip -0.05cm(\Cal E,p)\mapsto \bold{Bn}_0\, (\Cal E,p)& \text{on objects}\hskip -0.2cm\\
\hskip -0.05cm(F:\Cal E\to \Cal E')\mapsto \bold{Bn}_0\, (F):{\cases  \hskip -0.05cm(X,E)\mapsto (X,F(E)) & (X,E)\in Ob\, (\bold{Bn}_0\, (\Cal E,p)) \\
\hskip -0.05cm(f,\Phi )\mapsto (f,F_{E,E'}\circ \Phi )  & (f,\Phi )\in Ar\, (\bold{Bn}_0\, (\Cal E,p))\endcases }
& \text{on arrows}\hskip -0.2cm\endcases $
\vskip 0.2cm
and \ $\bold{Bn}\, (-):A\bold{Man}\text{-}\bold{CAT}\to \bold{Bn}\hookrightarrow 1\text{-}\bold{CAT}$,  such that \, $\bold{Bn}\, (-)=\bold{Man}(\bold{Bn}_0(-))$ (see 6, remarks)
(i.e., to each fibrewise functor there corresponds an actual functor on fibre bundles).
\endproclaim 
\demo{Proof} The given assignment for $\bold{Bn}_0\, (-)$ is obviously functorial. If $\Cal A:=\biggl \{\vcenter{\xymatrix{U_i\times p(E_i) \ar[d] \\ U_i}}\biggr \}_{i\in I}$ 
is an $\Cal E$-atlas for $\vcenter{\xymatrix{X \ar[d] \\ B}}\in Ob\, (\bold{Man}^{\to })$ then $\Cal A':=\biggl \{\vcenter{\xymatrix{U_i\times p'(F(E_i)) \ar[d] \\ U_i}}\biggr \}_{i\in I}$
is an $\Cal E'$-atlas for $\vcenter{\xymatrix{X' \ar[d] \\ B}}:=\biggl |\bold{Bn}\, (F)\biggl (\vcenter{\xymatrix{X \ar[d] \\ B}},\ \Cal A\biggr )\biggr |\in Ob\, (\bold{Man}^{\to })$ 
(i.e., essentially, $\Cal A'$ is a compatible family of arrows, if $\Cal A$ is compatible, and can be glued to an arrow $\vcenter{\xymatrix{X' \ar[d] \\ B}}$). 
\, So, $\bold{Bn}_0\, (F)$ and $\bold{Bn}\, (F)$ preserve atlases. \hfill $\square $
\enddemo 
\vskip 0.2cm
{\bf Remark.} Similarly, there can be defined fibrewise functors of more than one variables 
(e.g., $\bold{Bn}\, (\Cal E,p)\times \bold{Bn}\, (\Cal E',p')\to \bold{Bn}\, (\Cal E'',p'')$ induced by $F:\Cal E\times \Cal E'\to \Cal E''$, an $A\bold{Man}$-functor).
In this way usual fibrewise operations like $\times $, \, $\oplus $, \, $\otimes $, \, etc., are introduced.   \hfill $\square $

\head {\bf 6. Stacks and construction of general manifolds} \endhead

\vskip 0.1cm
Stacks give an example of relative higher order Category Theory. $n$-categories form an\linebreak 
$(n+1)$-category, so that (forgetting set-theoretical difficulties) 
$\text{Hom}_{(n+1)\text{-}\bold{CAT}}(\bold{C}, n\text{-}\bold{CAT})$ is an $n$-category.  

\vskip 0.2cm
{\bf Definition 6.1.} Let $(\bold{B},\tau )$ be a site ($\bold{B}$ is $1$-category), $F:\bold{B}^{op} \to n\text{-}\bold{CAT}$ be a (weak) functor.
\item{$\bullet $} For a sieve $i_s:s\hookrightarrow \bold{B}(-,B)$ \ $(n-1)$-category $\bold{Desc}(s,F):=\text{Hom}_{(n+1)\text{-}\bold{CAT}}(s,F)$ is 
called {\bf descent data} for functor $F$ and sieve $s$.
\item{$\bullet $} $\forall B\in Ob\, (\bold{B})$ and $\forall \text{ sieve }i_s:s\hookrightarrow \bold{B}(-,B)$ there is an induced ({\bf restriction}) 
functor $i_s^*:\text{Hom}_{(n+1)\text{-}\bold{CAT}}(\bold{B}(-,B),F)\to \text{Hom}_{(n+1)\text{-}\bold{CAT}}(s,F)$. 
If $i^*_s$ is full and faithfull ($\forall B \text{ and }\forall s$) then $F$ is called {\bf prestack}. If, moreover, it is an equivalence 
$F$ is called {\bf stack}, i.e.,  
$F$ is {\bf prestack} \ iff \ $\forall B, s$ \ $\text{Hom}_{(n+1)\text{-}\bold{CAT}}(s,F) @<i^*_s<\text{full, faith.}< \text{Hom}_{(n+1)\text{-}\bold{CAT}}(\bold{B}(-,B),F)@>\sim >\text{Yoneda}>F(B)$
and {\bf stack} iff $i^*_s$ is an equivalence.  \hfill $\square $

\vskip 0.1cm
For $n=1$ there is another definition of stack via matching families \cite{Moe, Vis}. \vskip 0.1cm 
Denote by $\text{\tt PreSt}(\bold{B}^{op},n\text{-}\bold{CAT}), \, \text{\tt St}(\bold{B}^{op},n\text{-}\bold{CAT})\hookrightarrow \text{Hom}_{(n+1)\text{-}\bold{CAT}}(\bold{B}^{op}, n\text{-}\bold{CAT})$ 
full subcategories of prestacks and staks respectively. 

\vskip 0.1cm
\proclaim{\bf Proposition 6.1} Both inclusions  \newline 
\vskip -0.3cm \hskip 1cm  $\xymatrix{\text{\tt St}(\bold{B}^{op},1\text{-}\bold{CAT}) \ar@<-1.1ex>@{^{(}->}[r]_-{i_0} & \text{\tt PreSt}(\bold{B}^{op},1\text{-}\bold{CAT}) \ar@<-1.1ex>@{^{(}->}[r]_-{i_1} \ar@{-->}@/_/[l]^-{\perp }_-{\Phi } & \text{Hom}_{\, 2\text{-}\bold{CAT}}(\bold{B}^{op}, 1\text{-}\bold{CAT}) \ar@{-->}@/_/[l]^-{\perp }_-{\Psi } }$
\item{}have left adjoints.
\endproclaim 
\demo{Proof} See \cite{Moe, Vis}.   \hfill $\square $
\enddemo

\vskip 0.35cm
\centerline {\bf Construction of manifolds of type $\bold{E}$ over site $(\bold{B}, \tau )$}

\item{1.} Factor $(1)$-functor $F:\bold{E}\to (\bold{B},\tau )$ through a free fibration (see {\bf Lemma 1.1}) $\vcenter{\xymatrix{\bold{E} \ar[r]^-{i} \ar[dr]_-{F} & 1/F \ar[d]^-{\text{dom}} \\
 & \bold{B}
}}$
\vskip -0.35cm
\item{2.} For fibration $\vcenter{\xymatrix{1/F \ar[d]_-{dom}\\
\bold{B}}}$ form a corresponding (weak) functor $\hat{F}:\bold{B}^{op}\to 1\text{-}\bold{CAT}$ and complete it to a stack
$(\Phi \circ \Psi )\hat{F}:\bold{B}^{op}\to 1\text{-}\bold{CAT}$ with respect to topology $\tau $.
\vskip 0.1cm
\item{3.} Get back (by Grothendieck construction) from stack $(\Phi \circ \Psi )\hat{F}:\bold{B}^{op}\to 1\text{-}\bold{CAT}$
to a fibration $\vcenter{\xymatrix{\tilde{\bold{E}} \ar[d]_-{p} & {\, 1/F} \ar@{_{(}->}[l] \ar[dl]^-{dom}  \\ \bold{B} & }}$
\vskip 0.2cm
\item{4.} Choose a  (correct) class of arrows $\Cal M$ in $\bold{B}$ representing 'embeddings of simple pieces into manifolds'.  
\vskip 0.1cm
\item{5.} Take a full subcategory $\bold{E}\text{-}\bold{Man}\hookrightarrow \tilde{\bold{E}}$ consisting of all objects 
$(B,\Cal E)\in Ob\, (\tilde{\bold{E}})$ such that $\exists \text{ a sieve } s\hookrightarrow \bold{B}(-,B)$ (depending on $B$) and 
$\forall f\in s$ $\bold{Cart}_f(\Cal E)=((\Phi \circ \Psi )\hat{F}(f))(\Cal E):(df\to F(E))\in \Cal M$ for some $E\in Ob\, (\bold{E})$. 
Then \ $\vcenter{\xymatrix{{\bold{E}\text{-}\bold{Man}\, } \ar@{^{(}->}@<-0.5ex>[r] \ar[dr]_-{p_F} & \tilde{\bold{E}} \ar[d]_-{p} & {\, 1/F} \ar@{_{(}->}@<0.5ex>[l] \ar[dl]^-{dom} \\
 & \bold{B} & }}$ \ is the required category of manifolds of type $\bold{E}$ over base site $(\bold{B},\tau )$. \hfill $\square $
 
\vskip 0.35cm
{\bf Remarks.} 
\vskip 0.1cm
\item{$\bullet $}  Depending on the choice of class $\Cal M$ categories $\bold{E}\text{-}\bold{Man}$ will be different 
(so, $\Cal M$ is an additional parameter). For cases of usual manifolds (smooth real or complex) $\Cal M$ is always class of topological 
embeddings of open subspaces. 
\item{$\bullet $} An object $(B,\Cal E)$ in $\bold{E}\text{-}\bold{Man}$ consists of a base object $B$ and an 'atlas' $\Cal E$, where $\Cal E$
is a class of compatible charts $(U\to F(E))\in \Cal M$, $U\in Ob\, (\text{Im}F)$, $E\in \bold{E}$. All arrows are represented by 
vertical arrows for the stack completion of $1/F$.
\item{$\bullet $} $\text{Im}(p_F)\supset \text{Im}(F)$.  
\vskip -0.2cm
\item{$\bullet $} The resulting category of manifolds $\vcenter{\xymatrix{\bold{E}\text{-}\bold{Man} \ar[d]^-{p_F} \\ \bold{B} }}$ 
is not usually a fibration.
\item{$\bullet $} Denote by $\bold{Man}_0\hookrightarrow 1\text{-}\bold{CAT}$ a category consisting of subcategory of $\bold{E}\text{-}\bold{Man}$ of trivial 
manifolds of type $\bold{E}$ for each type $\bold{E}$ and functors 'mapping $\bold{E}$-atlases to $\bold{E'}$-atlases'. Respectively, by 
$\bold{Man}\hookrightarrow 1\text{-}\bold{CAT}$ a category consisting of $\bold{E}\text{-}\bold{Man}$ for each type $\bold{E}$ and functors 
'mapping $\bold{E}$-atlases to $\bold{E'}$-atlases'. Then $\exists $ 'manifoldification' functor $\bold{Man}(-):\bold{Man}_0\to \bold{Man}$.
Inclusion functor $Ob\, (\bold{Man}_0)\ni \bold{E}\text{-}\bold{Man}_0\hookrightarrow \bold{E}\text{-}\bold{Man}\in Ob\, (\bold{Man})$ is natural in $\bold{E}$ (and itself is preserving atlases).  \hfill $\square $

\vskip 0.0cm
\centerline{{\bf Example} (single manifold)}

\vskip 0.1cm
Let $\bold{E}=$category consisting of $k^n, n=0,1,...$, ($k=\Bbb R$ or $\Bbb C$) and all smooth (analytic) local isomorphisms, 
$\Cal M$ be embeddings of open subspaces, 
$\bold{B}=$union of category of open subsets of a space $X$ with inclusion arrows, category $\bold{E}$, and arrows from $\Cal M$ with codomains in $\bold{E}$, 
$F:\bold{E}\to \bold{B}$ be the inclusion functor. Assignment 
$X\subset U\mapsto \{f:U\to k^n\, |\, n=0,1,..., f\in \Cal M\}$ gives a prestack on $X$. It is a nontrivial stack iff $X$ is a manifold.
$(X,\Cal E)\in Ob\, (\bold{E}\text{-}\bold{Man})$ iff $\Cal E$ is an atlas on $X$.

\head {\bf 7. Lifting problem for a group action} \endhead
Let $\bold{Grp}$ be a category of groups, $(-)_{\bold{Grp}}:1\text{-}\bold{CAT}\to 1\text{-}\bold{CAT}$ be a functor which assigns 
(in an obvious way) to each category without group action a category with groups actions, namely,
\vskip 0.2cm
$\cases \Cal C\mapsto \bold{Grp}\text{-}\Cal C & \ \, \Cal C\in Ob\, (1\text{-}\bold{CAT})\\
(F:\Cal C\to \Cal C')\mapsto (\bold{Grp}\text{-}F:\bold{Grp}\text{-}\Cal C\to \bold{Grp}\text{-}\Cal C') & \ F\in Ar\, (1\text{-}\bold{CAT})\endcases $ \hskip 1cm where:
\vskip 0.2cm
\hskip -0.39cm$\bold{Grp}\text{-}\Cal C$ consists of triples $(G,C,\rho )$ ($G\in Ob\, (\bold{Grp})$, $C\in Ob\, (\Cal C)$, $\rho :G\to \bold{Aut}(C)$ is a group homomorphism) as objects, 
and pairs $(\sigma :G\to G', f:C\to C'):(G,C,\rho )\to (G',C',\rho ')$ (s.t. $\forall g\in G \ \ \rho '(\sigma (g))\circ f=f\circ \rho (g)$\, ) as arrows, 
\vskip 0.1cm
$\bold{Grp}\text{-}F:\cases  (G,C,\rho )\mapsto (G,F(C),F_{C,C}\circ \rho )    &  \ \ (G,C,\rho )\in Ob\, (\bold{Grp}\text{-}\Cal C)  \\
(\sigma , f)\mapsto (\sigma ,F(f))   & \ \ \ (\sigma , f)\in Ar\, (\bold{Grp}\text{-}\Cal C)    \endcases $ \newline 
\vskip 0.1cm
\hskip -0.38cm[$(\sigma ,F(f))$ is an (equivariant) arrow in $\bold{Grp}\text{-}\Cal C'$ because $F(\rho '(\sigma (g)))\circ F(f)=F(f)\circ F(\rho (g) \ \ \forall g\in G$]

\proclaim{\bf Proposition 7.1} If $\vcenter{\xymatrix{\bold{E} \ar[d]^-{p} \\
\bold{B}}}$ is a {\bf structure over} $\bold{B}$ (i.e., all isomorphisms of type $(B'@>\sim >> p(E))\in Ar\, \bold{B}$ can be lifted to 
isomorphisms $(E'@>\sim >>E)\in Ar\, \bold{E}$) then $\vcenter{\xymatrix{\bold{Grp}\text{-}\bold{E} \ar[d]^-{\bold{Grp}\text{-}p} \\
\bold{Grp}\text{-}\bold{B}}}$ is a structure over $\bold{Grp}\text{-}\bold{B}$.
\endproclaim 
\demo{Proof} If $(\varphi ,f):(G',B',\rho ')@>\sim >>(G,p(E),p\circ \rho )$ is an iso then $\exists $  $\dbinom {E'}{B'}@>\sim >\hat f>\dbinom {E}{B}$, iso, 
because $\vcenter{\xymatrix{\bold{E} \ar[d]_-{p} \\
\bold{B}}}$ is a structure over $\bold{B}$. Regard the diagram (of group homomorphisms)
\vskip 0.1cm
\hskip 0.8cm$\vcenter{\xymatrix{\bold{Aut}_{\bold{E}}(E') \ar[dd]^-{ \ \ \ (**)}_-{p} &  &  & \bold{Aut}_{\bold{E}}(E) \ar[lll]_-{\hat f^{-1}\circ _{\bold{E}}-\circ _{\bold{E}}\hat f} \ar[dd]^-{p}\\
  & G'\ar[r]^-{\varphi }_-{\sim } \ar@{-->}[ul]_-{\rho ''} \ar[dl]^<(0.25){\rho '\hskip 1.25cm(*)} & G \ar[ur]^-{\rho } \ar[dr]_-{p\circ \rho } & \\
\bold{Aut}_{\bold{B}}(B') & & & \bold{Aut}_{\bold{B}}(p(E)) \ar[lll]^-{f^{-1}\circ _{\bold{B}}-\circ _{\bold{B}}f} }}$
\hskip 0.8cm$\rho ''(g'):=\hat f^{-1}\circ _{\bold{E}}\rho (\varphi (g'))\circ _{\bold{E}}\hat f$
\vskip 0.5cm
$(*)$ commutes because $f\circ _{\bold{B}}\rho '(g')=(p\circ \rho)(\varphi (g'))\circ _{\bold{B}}f$ by equivariance condition.
\vskip 0.1cm
$(**)$ commutes because $p(\rho ''(g'))=f^{-1}\circ _{\bold{B}}p(\rho (\varphi (g')))\circ _{\bold{B}}f=\rho '(g')$.
\vskip 0.1cm
So, $\exists $  iso $\dbinom {(G',E',\rho '')}{(G',B',\rho ')}@>{\binom {(\varphi ,\hat f)}{(\varphi , f)}}>\sim >\dbinom {G,E, \rho }{G,B,p\circ \rho }$, i.e.
$\vcenter{\xymatrix{\bold{Grp}\text{-}\bold{E} \ar[d]_-{\bold{Grp}\text{-}p} \\
\bold{Grp}\text{-}\bold{B}}}$ is a structure over $\bold{Grp}\text{-}\bold{B}$.   \hfill $\square $
\enddemo 

\vskip 0.0cm
There is a commutative diagram in $1\text{-}\bold{CAT}$ \ $\vcenter{\xymatrix{\bold{E} \ar[d]_-{p} & \bold{Grp}\text{-}\bold{E} \ar[l] \ar[d]^-{\bold{Grp}\text{-}p} \\
\bold{B} & \bold{Grp}\text{-}\bold{B} \ar[l] }}$ \ (where horizontal arrows forget group actions). So, there exists a forgetful 
fiber functor $\bold{E}_B@<<<\bold{Grp}\text{-}\bold{E}_{(G,B,\rho )}$.

\vskip 0.2cm
{\bf Definition 7.1.} For a given $G$-action $(G,B,\rho )\in Ob\, (\bold{Grp}\text{-}\bold{B}_B)$, an object $E\in Ob\, (\bold{E}_B)$
{\bf admits lifting of $G$-action} if $\exists $  $(G,E,\rho ')\in Ob\, (\bold{Grp}\text{-}\bold{E}_{(G,B,\rho )})$, i.e., 
$\vcenter{\xymatrix{E \ar@{|->}[d]_-{p} & (G,E,\rho ') \ar@{|->}[l] \ar@{|->}[d]^-{\bold{Grp}\text{-}p} \\
B & (G,B,\rho ) \ar@{|->}[l] }}$ (essentially, $\rho =p\circ \rho '$).   \hfill $\square $

\vskip 0.2cm
{\bf Lifting problem} for a $G$-action $\rho :G\to \bold{Aut}_{\bold{B}}(B)$ is equivalent to completion of the diagram of group homomorphisms with exact row
$\xymatrix{1 \ar[r] & \bold{Aut}_{\bold{E}_B}(E) \ar@{>->}[r] & \bold{Aut}_{\bold{E}}(E) \ar[r]^-{p} & \bold{Aut}_{\bold{B}}(B) \\
  &  &  G \ar@{-->}[u]_-{?}^-{\rho '} \ar[ur]_-{\rho } &  }$ 
where $\bold{Aut}_{\bold{E}_B}(E)$ are vertical automorphisms of $E$ over $B$.

\vskip 0.1cm
For single element $g\in \bold{Aut}_{\bold{B}}(B)$ there is a simple criterion of existance of $g'\in \bold{Aut}_{\bold{E}}(E)$ 
such that $p(g')=g$ (see \cite{Kom}). 

\proclaim{\bf Proposition 7.2} For a fibration $\vcenter{\xymatrix{\bold{E} \ar[d]_-{p} \\ 
\bold{B}}}$ (or structure over $\bold{B}$) $g\in \bold{Aut}_{\bold{B}}(B)$ can be lifted to $g'\in \bold{Aut}_{\bold{E}}(E)$ iff \,
$\bold{Cart}_g(E)@>\sim >>E$ ({\bf vertical} iso).
\endproclaim 
\demo{Proof} \ '$\Longleftrightarrow $' \hskip 2.5cm$\xymatrix{E \ar[dr]^<(0.6){\hskip -0.2cm\sim }_-{{\eightrm vertical}} \ar[drr]^<(0.55){g'}_<(0.55){\hskip 0.2cm\sim } & & \\
 & \bold{Cart}_g(E) \ar[r]^<(0.2){\sim }_<(0.2){\tilde g} & E\\
  & B \ar[r]^-{\sim }_-{g} & B}$  \vskip -0.399cm\hfill $\square $
\enddemo 

\vskip 0.0cm
\proclaim{\bf Proposition 7.3} If $\vcenter{\xymatrix{\bold{E} \ar[d] \\ \bold{B}}}$ is a structure on $\bold{B}$ (or 
(co)fibration with {\bf unique} (co)Cartesian lifting), and $(G,B,\rho )\in Ob\, (\bold{Grp}\text{-}\bold{B})$, then $\exists $ 
a representation \hskip 0.8cm$\xymatrix{G \ar[r]^-{\Cal R} \ar[dr]_-{\rho } & \bold{Aut}_{1\text{-}\bold{CAT}}(\bold{E}_B) \\
   &  \bold{Aut}_{\bold{B}}(B) \ar[u]}$ 
\vskip 0.0cm
where $\Cal R(g):\cases E\mapsto \bold{Cart}_{\rho (g)}(E)   &  \hskip 0.5cmE\in Ob\, (\bold{E}_B)     \\
\vcenter{\xymatrix{E \ar[d]_-{f}^-{\hskip 0.8cm\mapsto } & \bold{Cart}_{\rho (g)}(E) \ar[l]_-{\widetilde {\rho (g)}} \ar@{-->}[d]^-{\Cal R(g)(f)} \\
E' &  \bold{Cart}_{\rho (g)}(E') \ar[l]^-{\widetilde {\rho (g)}} }}    & \hskip 0.5cmf\in Ar\, (\bold{E}_B)    \endcases $
\endproclaim 
\vskip 0.05cm
\demo\nofrills{Proof\ \ \ } is straightforward.  \hfill $\square $
\enddemo

\vskip 0.1cm
{\bf Corollary.} If $E\in Ob\, (\bold{E}_B)$ is such that $\forall g\in G \ \bold{Cart}_{\rho (g)}(E)=E$ then Cartesian lifting
$\rho (g)\mapsto \widetilde {\rho (g)}$ lifts action $(G,B,\rho )\in Ob\, (\bold{Grp}\text{-}\bold{B})$ to the action 
$(G,E,\tilde {\rho })\in Ob\, (\bold{Grp}\text{-}\bold{E}_{(G,B,\rho )})$. \hfill $\square $

\vskip 0.35cm
\centerline{{\bf Example} (Covering Space)}
A covering space is a (co)fibration $\vcenter{\xymatrix{E \ar[d]_-{p} \\
B}}$ over groupoid $B$ with unique (co)cartesian lifting in which all morphisms are (co)cartesian. Moreover, representation $\bold{Aut}(b)\to \bold{Aut}(E_b)$, 
$b\in Ob\, (B)$ (induced by (co)cartesian ligting) is transitive on objects of $E_b$.

\proclaim{\bf Proposition 7.4} For a covering space $\vcenter{\xymatrix{E \ar[d]_-{p} \\
B}}$ over connected groupoid \, $B$
\item{$\bullet $} $\bold{Aut}\biggl (\vcenter{\xymatrix{E \ar[d]_-{p} \\
B}}\biggr )\simeq \bold{Aut}(E_b)$ \, (where \, $g\in \bold{Aut}(E_b)$ \, iff \, $g\circ f^*=f^*\circ g$, \, $f^*\equiv co\bold{Cart}_f$, \, $\forall f\in \bold{Aut}(b)$),
\item{$\bullet $} $\bold{Aut}(E_b)\simeq W(\bold{Stab}(e))\simeq N(\bold{Stab}(e))/\bold{Stab}(e)$ \, (where \, $\bold{Stab}(e)\hookrightarrow \bold{Aut}(b)$ \, is the stabilizer of an object \, $e\in Ob\, (E_b)$, \, $N(\bold{Stab}(e))$, $W(\bold{Stab}(e))$ \, are its normalizer and Weil group respectively).
\endproclaim 
\demo{Proof} 
\item{$\bullet $} An automorphism \, $g$ \, of covering space \, $p$ \, is given by family of fiberwise functors \, $g_b, \ b\in Ob\, (B)$, \, such that 
\, $f^*\circ g_b=g_{b'}\circ f^*$, \, $f^*\equiv co\bold{Cart}_f$, \, $\forall \, (f:b\to b')\in Ar\, (B)$. Take $g_b\in \bold{Aut}(E_b)$
and define \, $g_{b'}:=h^*\circ g_b\circ (h^*)^{-1}$ \, for some $h:b\to b'$. Then \, $g_{b'}$ \, is well-defined (if \, $h_1:b\to b'$ is 
another morphism then \, $h^*\circ g_b\circ (h^*)^{-1}=h_1^*\circ g_b\circ (h_1^*)^{-1}$ \, since $(h_1^{-1}\circ h)^*\circ g_b=(h_1^*)^{-1}\circ h^*\circ g_b=g_b\circ (h_1^*)^{-1}\circ h^*=g_b\circ (h_1^{-1}\circ h)^*$, \, $h_1^{-1}\circ h\in \bold{Aut}(b)$), 
and it is an automorphism of covering space \, $p$ \, (if \, $f:b'\to b''$ \, then \, $f^*\circ g_{b'}=g_{b''}\circ f^*$ \, since \, 
$f^*\circ h^*\circ g_b=g_{b''}\circ f^*\circ h^*$,  \, $f\circ h:b\to b''$). 
\item{$\bullet $} See \cite{May}.  \hfill $\square $
\enddemo

\vskip 0.3cm
\subhead 7.1. Lifting of a groupoid action for a sheaf\endsubhead

\vskip 0.2cm
{\bf Definition 7.1.1.} Let $(\bold{Top},\tau _0)$ be a site for all open coverings on topological spaces. 
\item{$\bullet $} $\bold{Set}$-valued presheaf $P:\bold{Top}^{op}\to \bold{Set}$ is a {\bf sheaf} iff $\forall $ sieve $S\hookrightarrow \bold{B}(-,B)$ and 
$\forall $ natural transformation $f:S\to P$ $\exists !$ $\hat f:\bold{B}(-,B)\to P$ such that \hskip 0.3cm$\xymatrix{{S\ } \ar@{^{(}->}[r] \ar[dr]_-{\forall f} & \bold{B}(-,B) \ar@{-->}[d]^-{\exists !\hat f}\\
  & P }$
\item{$\bullet $} $\bold{Cat}$-valued presheaf $P:\bold{Top}^{op}\to \bold{Cat}$ is a {\bf sheaf} iff its object and morphism parts 
are sheaves, i.e., \ $\bold{Top}^{op}@>P>>\bold{Cat}@>Ob>>\bold{Set}$ \ and \ $\bold{Top}^{op}@>P>>\bold{Cat}@>Mor>>\bold{Set}$ \ are
$\bold{Set}$-valued sheaves.
\item{$\bullet $} For presheaf $P:\bold{Top}^{op}\to \bold{Cat}$, space $X\in Ob\, (\bold{Top})$ and sieve $S\hookrightarrow \bold{Top}(-,X)$ matching family of objects 
$\tilde E:S\to Ob\circ P$ (nat. trans.) (or matching family of arrows $\tilde f:S\to Mor\circ P$ (nat. trans.)) {\bf has a germ} $\bold{germ}_x(\tilde E)$ (respectively, $\bold{germ}_x(\tilde f)$) at point $x\in X$ iff
$\exists $ $Colim_{s\in S,\, Im(s)\ni x}(\tilde E(s))=:\bold{germ}_x(\tilde E)$ (respectively, 
$Colim_{s\in S,\, Im(s)\ni x}(\tilde f(s))=:\bold{germ}_x(\tilde f)$) (if germ exists it is unique up to iso and does not depend on the 
choice of sieve $S$).
\item{$\bullet $} {\bf Etale space} is $E:=\coprod \limits _{x\in X}\bold{germ}_x(\tilde E)$ (respectively, $f:=\coprod \limits _{x\in X}\bold{germ}_x(\tilde f)$) 
(depending on two variables: 'point' $x\in X$ and 'matching family' $\tilde E$ or $\tilde f$) with topology generated by basic open
sets $\bigl (U,\{\bold{germ}_x(\tilde E)\, |\, x\in U\}\bigr )$ \, (or, $\bigl (U,\{\bold{germ}_x(\tilde f)\, |\, x\in U\}\bigr )$),
\, $U$ is open in $X$. There is a natural continuous projection $p:E\to X:(x,\bold{germ}_x(\tilde E))\mapsto x$ (respectively, $p:f\to X:(x,\bold{germ}_x(\tilde f))\mapsto x$)
which is a local homeomorphism.    \hfill $\square $

\vskip 0.2cm
\proclaim{\bf Lemma 7.1.1} Every fibration is a cofibration with respect to iso's (every cofibration is a fibration with 
respect to iso's).
\endproclaim 
\demo{Proof} Let $\vcenter{\xymatrix{\bold{E} \ar[d]_-{p}\\ 
\bold{B} }}$ be a fibration, and $p(E)@>f>> B'$ be an iso in $\bold{B}$.  Then $\tilde f:=\bigl (\widetilde {f^{-1}}\bigr )^{-1}:E\to E'$ (where $\widetilde { \ \ }$ on the right is a cartesian lifting) is a 
cocartesian lifting of $f$ (obvious).  \hfill $\square $
\enddemo 

{\bf Corollary.} For a (co)fibration $\vcenter{\xymatrix{\bold{E} \ar[d]_-{p} \\
\bold{B}}}$ for each iso $(f:B@>\sim >>B')\in Ar\, \bold{B}$ inverse image $\bold{Cart}_f:\bold{E}_{B'}\to \bold{E}_{B}:s_{B'}\mapsto f^*(s_{B'})$ 
and direct image $\bold{Cart}_f:\bold{E}_{B}\to \bold{E}_{B'}:s_{B}\mapsto f_*(s_{B})$ \, (where $s_B$ is a 'section' (object or morphism) over $B$) 
\, exist.   \hfill $\square $

\vskip 0.2cm
{\bf Definition 7.1.2.} 
\item{$\bullet $} For a space $X\in Ob\, (\bold{Top})$ {\bf groupoid of local homeomorphisms} of $X$ is a subcategory $\bold{Gr}_X\hookrightarrow \bold{Top}$ 
such that $\cases  Ob\, (\bold{Gr}_X)  & \text{are open subsets in }X \\
Ar\, (\bold{Gr}_X)     & \text{are iso's in }\bold{Top} \text{ \ (between open subsets in }X\text{)}\endcases $
\vskip 0.1cm
(Nonfull) subcategory $\bold{Gr}_{X,x}\hookrightarrow \bold{Gr}_X$ with objects $U\ni x$ and morphisms $f:U\to V, f(x)=x$,
\vskip 0.0cm
is called
{\bf groupoid of local homeomorphisms} of $X$ {\bf with fixed point} $x\in X$.
\item{$\bullet $} $X$ is {\bf transitive} with respect to $\bold{Gr}_X$ if $\forall \, x,y\in X$ $\exists \, U,V\in Ob\, (\bold{Gr}_X),\, (f:U\to V)\in Ar\, (\bold{Gr}_X)$ 
such that $U\ni x, V\ni y,\, f(x)=y$. 
\item{$\bullet $} For a (co)fibration $\vcenter{\xymatrix{\bold{E} \ar[d]_-{p} \\
\bold{Top}}}$ with unique (co)cartesian ligting and space $X\in Ob\, (\bold{Top})$ two actions of $\bold{Gr}_X$ on local sections 
over $X$ are defined:
\vskip 0.1cm
{\bf left action} $\forall f\in \bold{Gr}_X(U,V)$ \ $f^*\equiv \bold{Cart}_f:\bold{E}_V\to \bold{E}_U:s_V\mapsto f^*s_V$
\vskip 0.1cm 
{\bf right action} $\forall f\in \bold{Gr}_X(U,V)$ \ $f_*\equiv co\bold{Cart}_f:\bold{E}_U\to \bold{E}_V:s_U\mapsto f_*s_U$
\vskip 0.1cm
(where $s_V,s_U$ are objests or morphisms).
\item{$\bullet $} To each of actions $f^*,f_*$ (on local sections of $\bold{E}_X$) there correspond respectively left and right
actions of $\bold{Gr}_{X,x}$ on $\{\, \bold{germ}_x(\tilde s)\ |\ \tilde s \text{ is a matching family of local sections of } \bold{E}_X\, \}$.
If $\goth s=\bold{germ}_x(s_U)$ is a germ at point $x$ presented by a local section $s_U$ (i.e., $\goth s=\underset {U\supset V\ni x}\to {Colim}(s_U\underset V\to |)$) 
then if \, $(f:W\to V)\in \bold{Gr}_{X,x}(W,V)$
\vskip 0.15cm
{\bf left action} \ $f^*\goth s:=\bold{germ}_x(\bigl (f\underset {f^{-1}(U\cap V)}\to {|}\bigr )^*\bigl (s_U\underset {U\cap V}\to {|}\bigr ))$ 
\vskip 0.1cm
{\bf right action} \ $f_*\goth s:=\bold{germ}_x(\bigl (f\underset {U\cap W}\to {|}\bigr )_*\bigl (s_U\underset {U\cap W}\to {|}\bigr ))$. 
\item{$\bullet $} For a subgroupoid $\bold{G}\hookrightarrow \bold{Gr}_X$, \, $\dbinom sX\in \biggl (\vcenter{\xymatrix{\bold{E} \ar[d]_-{p} \\
\bold{Top}}}\biggr )$ is {\bf $\bold{G}$-invariant} 
if $\forall \, (f:U\to V)\in Ar\, (\bold{G})$ $f^*(s_V)=s_U$ (or $f_*(s_U)=s_V$) (where: $s_U:=s\underset U\to |:=(i:U\hookrightarrow X)^*s$,  
$\vcenter{\xymatrix{\bold{E} \ar[d]_-{p} \\
\bold{B}}}$ is a fibration with unique cartesian lifting (or a local structure with respect to inclusions of open sets), $s$ is a section (object or morphism) over $X$).
In other words, $\bold{G}$-invariant sections {\bf admit lifting} of groupoid $\bold{G}$. 
\item{$\bullet $} $\bold{germ}_{x,y}(f)$ of a map $(f:U\to V)\in Ar\, (\bold{Gr}_X)$, such that $f(x)=y$, is an equivalence class of 
maps $\{\, g\in Ar\, (\bold{Gr}_X)\, |\, g(x)=y, \ \exists \text{ opens }W_x\ni x,W_y\ni y, \text{ such that } \exists \text{ the same restrictions }$\newline 
$f\underset {W_xW_y}\to |=g\underset {W_x,W_y}\to |\in Ar\, (\bold{Gr}_X) \, \}$. Assume, \, $\goth s_x=\bold{germ}_x(s_{U_1}), \, \goth s'_y=\bold{germ}_y(s'_{V_1})$. Then \newline 
{\bf left action} \ $(\bold{germ}_{x,y}(f))^*\goth s'_y:=\bold{germ}_x(\bigl (f\underset {f^{-1}(V\cap V_1)}\to |\bigr )^*\bigl (s'\underset {V\cap V_1}\to |\bigr ))$\newline 
{\bf right action} \ $(\bold{germ}_{x,y}(f))_*\goth s_x:=\bold{germ}_y(\bigl (f\underset {(U\cap U_1)}\to |\bigr )_*\bigl (s\underset {U\cap U_1}\to |\bigr ))$.
\hfill $\square $

\proclaim{\bf Lemma 7.1.2} 
\item{$\bullet $} Let $X$ be a topological space, $\bold{G}\hookrightarrow \bold{Gr}_X$ be a subgroupoid, $\bold{G}_x\hookrightarrow \bold{G}$ be a subgroupoid of pointed maps with fixed point $x\in X$. 
Then $\forall x,y\in X$ and $\forall f\in Ar\, (\bold{G})$, s.t. $f(x)=y$, $\bold{germ}_{x,x}(\bold{G}_x)=\bold{germ}_{y,x}(f^{-1})\cdot \bold{germ}_{y,y}(\bold{G}_y)\cdot \bold{germ}_{x,y}(f)$ (for certain unique composite $\cdot $ of germs of maps).
\item{$\bullet $} If $\goth s_x=\bold{germ}_x(s_U)\in S_x\subset S$ is a point of etale space $\vcenter{\xymatrix{S \ar[d] \\
X}}$ (corresponding to objects or morphisms over $X$ for a fibration $\vcenter{\xymatrix{\bold{E} \ar[d]_-{p} \\
\bold{Top}}}$\hskip -0.2cm) \ then $\goth s_x$ is $\bold{G}_x$-invariant iff it is $\bold{germ}_{x,x}(\bold{G}_x)$-invariant. 
\item{$\bullet $} If $\bold{G}$ is transitive on $X$ and $\goth s_x$ is $\bold{G}_x$-invariant then $\forall f,g\in Ar\, (\bold{G})$, s.t. $f(x)=y,\ g(x)=y$, 
there is a unique induced germ at point $y$ $(\bold{germ}_{x,y}(f))_*\goth s_x=(\bold{germ}_{x,y}(g))_*\goth s_x$ and this germ is $\bold{G}_y$-invariant 
(or, respectively, $(\bold{germ}_{y,x}(f^{-1}))^*\goth s_x=(\bold{germ}_{y,x}(g^{-1}))^*\goth s_x$ is a unique $\bold{G}_y$-invariant germ at point $y$), 
i.e. $\goth s_x$ can be distributed in a unique way over all $X$ (to give rise a section $s:X\to S$ of etale space $\vcenter{\xymatrix{S \ar[d] \\
X}}$ consisting of invariant germs at each point). \hfill $\square $
\endproclaim

\proclaim{\bf Proposition 7.1.1} For sheaf $P:\bold{Top}^{op}\to \bold{CAT}$, space $X\in Ob\, (\bold{Top})$, and transitive 
groupoid $\bold{G}\hookrightarrow \bold{Gr}_X$ $\bold{G}$-invariant sections over $X$ are in bijective correspondence with a subset
of $\bold{G}_x$-invariant germs (of local sections) for a fixed point $x\in X$.
\endproclaim
\demo{Proof} To each $\bold{G}$-invariant section over $X$ there corresponds $\bold{G}_x$-invariant germ of this section at point $x$.
Conversely, by lemma 7.1.2, each $\bold{G}_x$-invariant germ generates a section of the corresponding etale space. When this section
is continuous there is a global section over $X$ of sheaf $P$ (which is locally invariant). \hfill $\square $
\enddemo

{\bf Remark.} B.P. Komrakov \cite{Kom} asserts (without a proof) that the above bijective correspondence is with the whole set
of $\bold{G}_x$-invariant germs. But, without additional assumptions it is not clear why the corresponding section of invariant germs 
is continuous and the sheaf section is invariant.  \hfill $\square $ 

\head {\bf 8. Equivalence, groups, actions} \endhead

Let $\Cal R$ be a category of sets with a given equivalence relation for each set. There are following functors
\item{$\bullet $} {\bf forgetful} $p:\Cal R\to \bold{Set}:(A,R)\mapsto A$
\item{$\bullet $} {\bf quotient} $Q:\Cal R\to \bold{Set}:(A,R)\mapsto A/R$
\item{$\bullet $} {\bf inclusion} $\Delta :\bold{Set}\to \Cal R:A\mapsto (A,\Delta _A)$, $\Delta _A:=\{(a,a)\, |\, a\in A\}$
\item{} such that $\xymatrix{\Cal R \ar@/^/[r]^-{Q}_-{\perp } & \bold{Set} \ar@/^/[l]^-{\Delta }}$, i.e. $\bold{Set}(Q(A,R),B)\underset {nat.iso}\to {@>\sim >>}\Cal R((A,R),\Delta (B)):f\mapsto f\circ \pi $, where $\xymatrix{A \ar@{->>}[r]^-{\pi } & A/R}$ 
is the canonical projection (so, quotient object $Q(A,R)$ represents functor $\Cal R((A,R),\Delta (-)):\bold{Set}\to \bold{Set}$).

\vskip 0.2cm

For arbitrary category $\bold{C}$ equivalence relation on objects is introduced as usual via hom-sets.

\vskip 0.2cm

{\bf Definition 8.1.} \vskip 0.1cm
\item{$\bullet $} A functor $R:\bold{C}^{op}\to \Cal R$ is called an {\bf equivalence relation on object} $C\in Ob\, \bold{C}$
\ iff \,  \hskip 0.35cm$\vcenter{\xymatrix{\bold{C}^{op} \ar[dr]_-{\bold{C}(-,C)} \ar[r]^-{R} & \Cal R \ar[d]^-{p} \\
 & \bold{Set}}}$ \hskip 0.2cm(i.e. usual equivalence relations are introduced on hom-sets $\bold{C}(C',C), \ C'\in Ob\, \bold{C}$ \, and 
 they are preserved under precomposition \, $-\circ f, \ f:C''\to C'$).   
\item{$\bullet $} Let $\bold{C}_{\Cal R}$ be a category such that
\item{} $Ob\, (\bold{C}_{\Cal R})$ are pairs $(C,R)$, \, $C\in Ob\, \bold{C}$, \, $R$ is an equivalence relation on $C$,            
\item{} $Ar\, (\bold{C}_{\Cal R})$ are maps $(C,R)@>(f,F)>>(C',R')$, where $(f:C\to C')\in Ar\, \bold{C}$ and $F:R\Rightarrow  R'$ is a
natural transformation of equivalence relations such that \hskip 0.35cm$\xymatrix{\bold{C}^{op} \ar@/^/[r]^-{R}_-{\Downarrow F} \ar@/_/[r]_-{R'} \ar@/^2pc/[rr]^-{\bold{C}(-,C)} \ar@/_2pc/[rr]_-{\bold{C}(-,C')} & \Cal R \ar[r]_-{p} & \bold{Set}}$
\item{} $p_F=\bold{C}(-,f)$ \ (it means $(C,R)@>(f,F)>>(C',R')$ is a morphism in $\bold{C}_{\Cal R}$ iff $f:C\to C'$ is an arrow in $\bold{C}$ and 
$f\circ -$ preserves equivalence relation, i.e. if $g_1\sim _{R}g_2$ then $f\circ g_1\sim _{R'}f\circ g_2$ for $g_1,g_2:X\to C$).  \hfill $\square $

\vskip 0.2cm
$\bold{C}_{\Cal R}$ is an analogue of $\Cal R$ for arbitrary category $\bold{C}$. 
Again, there are the following functors
\item{$\bullet $} {\bf forgetful} $p:\bold{C}_{\Cal R}\to \bold{C}:(C,R)\mapsto C$
\item{$\bullet $} {\bf inclusion} $\Delta :\bold{C}\to \bold{C}_{\Cal R}:C\mapsto (C,\, \Delta \circ \bold{C}(-,C))$, where 
$\Delta :\bold{Set}\to \Cal R:A\mapsto (A,\Delta _A)$, $\Delta _A:=\{(a,a)\, |\, a\in A\}$
\item{$\bullet $} {\bf quotient} $Q:\bold{C}_{\Cal R}\to \bold{C}:(C,R)\mapsto C/R$ which is a left adjoint to $\Delta :\bold{C}\to \bold{C}_{\Cal R}$,
i.e. $\xymatrix{\bold{C}_{\Cal R} \ar@/^/[r]^-{Q}_-{\perp } & \bold{C} \ar@/^/[l]^-{\Delta } }$ \, or \,  
$\bold{C}(Q(C,R),\, C')\underset {nat. iso}\to {@>\sim >>}\bold{C}_{\Cal R}((C,R),\, \Delta (C'))$ \, 
(quotient object $C/R:=Q(C,R)$ represents functor $\bold{C}_{\Cal R}((C,R),\Delta (-)):\bold{C}\to \bold{Set}$ which means  
that $\exists $ an arrow $\pi :(C,R)\to \Delta (Q(C,R))$ \, such that \, $\forall f:(C,R)\to \Delta (C')$ \, $\exists ! \hat f:Q(C,R)\to C'$ 
\, with \, $f=\Delta (\hat f)\circ \pi $, in other words, quotient map $\pi :C\to Q(C,R)$ is a common coequalizer of all 
equivalent arrows $f\sim _{R}g$ with arbitrary domain and codomain $C$ and, in particular, is always an epimorphism).    

\vskip 0.1cm
\item{} Quotient functor may not exist for the whole category $\bold{C}_{\Cal R}$, but there always exists a (maximal) full 
subcategory $\bold{C}_{\Cal RQ}\hookrightarrow \bold{C}_{\Cal R}$ for which $\xymatrix{\bold{C}_{\Cal RQ} \ar@/^/[r]^-{Q}_-{\perp } & \bold{C} \ar@/^/[l]^-{\Delta } }$
(indeed, $\bold{C}_{\Cal RQ}$ is always non empty since $Q\circ \Delta (C)=C$, i.e. $\Delta (C)\in Ob\, (\bold{C}_{\Cal RQ})$).

\vskip 0.2cm
If $\bold{C}$ is a concrete category with representable underlying functor $U:=\bold{C}(I,-)$ then to each equivalence relation
$R:\bold{C}^{op}\to \Cal R$ on object $C$ with quotient map $\pi :(C,R)\to \Delta (Q(C,R))$ there corresponds a usual equivalence relation on $\bold{C}(I,C)$
with quotient map $\pi \circ -:\bold{C}(I,C)\to \bold{C}\bigl (I,Q(C,R)\bigr )$, and, conversely, to usual equivalence relation on
$\bold{C}(I,C)$ with quotient map $\pi \circ -:\bold{C}(I,C)\to \bold{C}(I,C')$ there corresponds a maximal 'saturated' equivalence 
relation $R:\bold{C}^{op}\to \Cal R$ on object $C$ with quotient map $\pi :C\to C'\equiv Q(C,R)$ such that $f\sim _{R}g$ iff $\pi \circ f=\pi \circ g$. In general, equivalence relation
on hom-sets is weaker than usual one.

\vskip 0.2cm
Let $C\in Ob\, \bold{C}$, \, $\sigma :G\to \bold{Aut}_{\bold{C}}(C)$, \, then $G$ also acts on hom-sets $\bold{C}(C',C), \ C'\in Ob\, \bold{C}$, \, 
$G\times \bold{C}(C',C)\to \bold{C}(C',C):\cases  (g,f)\mapsto \sigma (g)\circ f  &  \text{left action} \\
(g,f)\mapsto \sigma (g^{-1})\circ f & \text{right action}     \endcases $\, , \ i.e. $\exists $ \, a functor 
\, $\Sigma :\bold{C}^{op}\to G\text{-}\bold{Set}$ \, such that \hskip 0.2cm$\vcenter{\xymatrix{\bold{C}^{op} \ar[r]^-{\Sigma } \ar[dr]_-{\bold{C}(-,C)} & G\text{-}\bold{Set} \ar[d]^-{p} \\
  & \bold{Set}}}$ \hskip 0.2cm(it means that all hom-sets $\bold{C}(C',C), \ C'\in Ob\, \bold{C}$, are regarded with the given $G$-action).

\vskip 0.2cm
There are functors\vskip 0.1cm
\item{$\bullet $} $r:G\text{-}\bold{Set}\to \Cal R:\cases  (X,G,\sigma )\mapsto (X,R_{\sigma })  &  \text{on objects} \\
((X,G,\sigma )@>f>>(X',G,\sigma '))\mapsto ((X,R_{\sigma })@>f>>(X',R_{\sigma '}))  & \text{on arrows}  \endcases $
\vskip 0.1cm
(where $R_{\sigma }$ is an equivalence relation on $X$ such that $(x,y)\in R_{\sigma }$ iff $\exists \, g\in G \ y=\sigma (g)x$)
\item{} $r$ is a functor over $\bold{Set}$, i.e. \hskip 0.5cm$\vcenter{\xymatrix{G\text{-}\bold{Set} \ar[r]^-{r} \ar[dr]_-{p} & \Cal R \ar[d]^-{p} \\
     &  \bold{Set}}}$
\item{$\bullet $} $r:G\text{-}\bold{C}\to \bold{C}_{\Cal R}:\cases  (C,G,\sigma )\mapsto (C,R_{\sigma })   &   \text{on objects}      \\
((C,G,\sigma )@>f>>(C',G,\sigma '))\mapsto ((C,R_{\sigma })@>f>>(C',R_{\sigma '}))  & \text{on arrows}   \endcases $
\vskip 0.1cm
(where $R_{\sigma }:=r\circ \Sigma $  \, is an equivalence relation on object $C$ corresponding to \, $\sigma $,  \ $\vcenter{\xymatrix{\bold{C}^{op} \ar[dr]_-{\bold{C}(-,C)} \ar[r]^-{R_{\sigma }} & \Cal R \ar[d]^-{p} \\
           &  \bold{Set} }}$)
\item{} $r$ is a functor over $\bold{C}$, i.e. \hskip 0.5cm$\vcenter{\xymatrix{G\text{-}\bold{C} \ar[r]^-{r} \ar[dr]_-{p} & \bold{C}_{\Cal R} \ar[d]^-{p} \\
     &  \bold{C}}}$

\item{}Let $G\text{-}\bold{C}_{Q}:=r^{-1}(\bold{C}_{\Cal RQ})$. Then $\exists $ a quotient functor $G\text{-}\bold{C}_Q@>r>>\bold{C}_{\Cal RQ}@>Q>>\bold{C}$.
Denote it again by $Q$, and $Q\circ r\, (C,G,\sigma )$ by $C/G$.

\vskip 0.2cm
For arbitrary functor $F:\bold{C}\to \bold{D}$ we have $G\text{-}F:G\text{-}\bold{C}\to G\text{-}\bold{D}$ such that 
$\vcenter{\xymatrix{\bold{C} \ar[r]^-{F} & \bold{D} \\
G\text{-}\bold{C} \ar[r]_-{G\text{-}F} \ar[u]^-{p} &   G\text{-}\bold{D} \ar[u]_-{p} }}$
but $F$ needs not preserve quotients, i.e. the diagram \hskip 0.35cm$\vcenter{\xymatrix{G\text{-}\bold{C}_{Q} \ar@{..>}[r]^-{G\text{-}F} \ar[d]_-{Q}^-{ \ \ \ \ \ \ \ \cong }  & G\text{-}\bold{D}_{Q} \ar[d]^-{Q} \\
\bold{C} \ar[r]_-{F} &   \bold{D} }}$\hskip 0.35cm
can be wrong (the dotted arrow may not exist and the natural isomorphism may not hold). If the above diagram holds (up to iso)
then \, $F:\bold{C}\to \bold{D}$ \, {\bf preserves quotients} (of category $G\text{-}\bold{C}$). In this case $F(C/G)\cong F(C)/G$. 
Quotient $C/G$ is called {\bf universal} \cite{Kom} if $\forall C'\in Ob\, \bold{C}$ $\vcenter{\xymatrix{C'\times C \ar[d]_-{\pi } \ar[dr]^-{1\times p} & \\
(C'\times C)/G \ar[r]^-{\sim } & C'\times (C/G)\hskip -0.2cm}}$
\vskip -0.5cm
\hskip -0.38cm($C'$ is with trivial $G$-action).

\proclaim{\bf Proposition 8.1} Let $\vcenter{\xymatrix{\bold{E} \ar[d]_-{p} \\
\bold{B}}}$ be a structure on $\bold{B}$, $p$ preserve quotients of category $G\text{-}\bold{E}$, and $(E,G,\sigma )\in Ob\, (G\text{-}\bold{E})$ be 
an object such that $E/G$ exists with $\pi :E\to E/G$, canonical projection, then \, $\dbinom {E/G}{p(E/G)}=(p(\pi ))_*\dbinom {E}{p(E)}$
\, is a direct image of \, $\dbinom {E}{p(E)}$.
\endproclaim
\demo{Proof} We need to prove that $\dbinom {E}{p(E)}@>{\dbinom {\pi }{p(\pi )}}>>\dbinom {E/G}{p(E/G)}$ is cocartesian. 
Take $\dbinom {u}{v}\hskip -0.05cm:\hskip -0.05cm\dbinom {E}{p(E)}\hskip -0.05cm\to \dbinom {E'}{p(E')}$ such that $v=k\circ p(\pi )$ for some $k:p(E/G)\to p(E')$, i.e. 
$\forall \, f\sim _{G}f':B\to p(E) \ \ \, v\circ f=v\circ f'$.
\vskip 0.0cmAssume, $h\sim _{G}h':E_1\to E$ then $p(h)\sim _{G}p(h'):p(E_1)\to p(E)$ (because, if $h'=\sigma (g)\circ h$ then $p(h')=p(\sigma (g))\circ p(h)$).
So, $v\circ p(h)=v\circ p(h')$, \, $p(u)\circ p(h)=p(u)\circ p(h')$, \, $u\circ h=u\circ h'$ ($p$ is faithful), i.e. $u$ coequalizes
all \, $\sim _{G}$-equivalent arrows (in $R_{\sigma }$). 
\vskip 0.0cmTherefore, $u=\hat u\circ \pi $ for a unique $\hat u:E/G\to E'$.
\vskip 0.0cm$\bigl (v=p(u)=p(\hat u)\circ p(\pi )=k\circ p(\pi )\bigr )$ $\Rightarrow $ $\bigl (p(\hat u)=k\bigr )$ by universality of $p(\pi )$.
\vskip 0.0cmFinally, $\exists ! \, \dbinom {\hat u}{k}:\dbinom {E/G}{p(E/G)}\to \dbinom {E'}{p(E')}$ such that $\dbinom {u}{v}=\dbinom {\hat u}{k}\circ \dbinom {\pi }{p(\pi )}$, i.e.
$\dbinom {\pi }{p(\pi )}$ is cocartesian.               \hfill $\square $
\enddemo

\vskip 0.2cm
\subhead 8.1. Group objects, subgroups, quotient objects \endsubhead
\vskip 0.2cm
{\bf Definition 8.1.1.} Let $\bold{C}$ be a category with binary products and terminal object $1$.
\item{$\bullet $} $G\in Ob\, \bold{C}$ is called a {\bf group object} if $\exists $ maps $m:G\times G\to G$, \, $e:1\to G$, \, $inv:G\to G$
such that the following group-like diagrams hold \vskip 0.0cm \hskip -0.1cm
$\xymatrix{G\times G\times G \ar[d]_-{m\times 1} \ar[r]^-{1\times m} & G\times G \ar[d]^-{m} \\
G\times G \ar[r]_-{m}   &  G  }$ \hskip 0.0cm
$\xymatrix{G\times 1 \ar[r]^-{1\times e} \ar[dr]_-{p_1} & G\times G \ar[d]^-{m} & 1\times G \ar[l]_-{e\times 1} \ar[dl]^-{p_2} \\
  &  G }$ \hskip 0.0cm
$\xymatrix{G\times G \ar[r]^-{1\times inv} & G\times G \ar[d]^-{m} & G\times G \ar[l]_-{inv\times 1}\\
G \ar[u]^-{\Delta } \ar[r]_-{e\circ !} & G & G \ar[u]_-{\Delta } \ar[l]^-{e\circ !}\hskip -0.38cm}$  
\item{$\bullet $} Subobject $\xymatrix{K \ar@{>->}[r] & G}$ of group object $G$ is called a {\bf subgroup} (object) if $\exists $
maps $m_K:K\times K\to K$, \ $e_K:1\to K$, \ $inv_K:K\to K$ such that 
$\xymatrix{K\times K \ar@{>->}[r] \ar@/^1.25pc/[rrr]^-{m_K} & G\times G \ar[r]_-{m} & G & K \ar@{>->}[l]}$
$\xymatrix{1 \ar[r]_-{e_K} \ar@/^1.25pc/[rr]^-{e} & K \ar@{>->}[r] & G }$ \hskip 0.7cm
$\xymatrix{K \ar@{>->}[r] \ar@/^1.25pc/[rrr]^-{inv_K} & G \ar[r]_-{inv} & G & K \ar@{>->}[l] }$
\vskip -0.5cm\item{$\bullet $} For two elements \, $f,g:1\to G$ \, {\bf multiplication} \, $f\cdot g:1\to G$ \ is \hskip 0.25cm$\xymatrix{1 \ar[r]_-{<f,g>} \ar@/^1.25pc/[rr]^-{f\cdot g} & G\times G \ar[r]_-{m} & G}$
\item{$\bullet $} {\bf Right shift} \, $R_g:G\to G$ \, (by element \, $g:1\to G$) \, is \hskip 0.25cm$\xymatrix{G \ar[r]_-{\sim } \ar@/^1.25pc/[rrr]^-{R_g} & G\times 1 \ar[r]_-{1\times g} & G\times G \ar[r]_-{m} & G}$ 
{\bf Left shift} \, $L_g:G\to G$ \, (by element \, $g:1\to G$) \, is \hskip 0.2cm$\xymatrix{G \ar[r]_-{\sim } \ar@/^1.25pc/[rrr]^-{L_g} & 1\times G \ar[r]_-{g\times 1} & G\times G \ar[r]_-{m} & G}$       \hfill $\square $

\vskip 0.3cm
\proclaim{\bf Proposition 8.1.1} For a group object $G\in Ob\, \bold{C}$
\item{$\bullet $} $\bold{C}(1,G)$ is a group,
\item{$\bullet $} $\exists $ (anti)representation $\bold{C}(1,G)\to \bold{Aut}_{\bold{C}}(G):g\mapsto R_g$ (by right shifts)
and representation $\bold{C}(1,G)\to \bold{Aut}_{\bold{C}}(G):g\mapsto L_g$ (by left shifts).
\endproclaim 
\demo{Proof} 
\item{$\bullet $} It follows immediately from group object axioms ($e:1\to G$ is the identity, $inv\circ g:1\to G$ is the inverse of $g:1\to G$).
\item{$\bullet $} $R_e=1_G:G\to G$ (obvious)\vskip 0.0cm
$R_f\circ R_g=R_{g\circ f}$ follows from the diagram\vskip 0.0cm
\hskip -0.1cm$\xymatrix{G \ar[r]^-{\sim } \ar[dr]_-{\sim } & G\times 1 \ar@{..>}[r]^-{1\times g} & G\times G \ar@{..>}[r]^-{m} & G \ar@{..>}[r]^-{\sim } & G\times 1 \ar[r]^-{1\times f} & G\times G \ar[r]^-{m} & G \\
 & G\times 1  \ar[rrrr]^-{1\times (g\cdot f)} \ar@{..>}[dl]_-{1\times <1,1>} \ar[drr]^(0.55){1\times <g,f>} & & & & G\times G \ar[ur]_-{m} & \\
G\times 1\times 1 \ar@{..>}[drrr]_-{(m\circ (1\times g))\times 1} \ar[rrr]^-{1\times g\times f} & & & G\times G\times G \ar[urr]^(0.39){1\times m} \ar[rrr]^-{m\times 1} & & & G\times G \ar[uu]_-{m} \\
 & & & G\times 1 \ar[urrr]_-{1\times f} & & & \hskip -0.2cm}$
 \item{}Two dotted paths $\xymatrix{G\times 1 \ar@{..>}[r]^-{1\times g} & G\times G \ar@{..>}[r]^-{m} & G \ar@{..>}[r]^-{p_1^{-1}}_-{\sim } & G\times 1}$ and 
 \item{}$\xymatrix{G\times 1 \ar@{..>}[r]^-{1\times <1,1>} & G\times 1\times 1 \ar@{..>}[rr]^-{(m\circ (1\times g))\times 1} & & G\times 1 }$ are equal 
 since their composites with projections $p_1:G\times 1@>\sim >> G$ and $p_2\equiv \, !:G\times 1@>>>1$ are equal, indeed, \hskip -0.1cm
 $\cases \hskip -0.05cmp_1\circ (p_1^{-1}\circ m\circ (1\times g))=m\circ (1\times g)  & \hskip -0.5cm  \\
 \hskip -0.05cmp_2\circ (p_1^{-1}\circ m\circ (1\times g))=\, !  &  \hskip -0.5cm \endcases $ 
 and \hskip -0.05cm$\cases  \hskip -0.05cmp_1\circ (((m\circ (1\times g))\times 1)\circ (1\times <1,1>))=p_1\circ (((m\circ (1\times g))\times 1)\circ (<1,!>\times 1))=*   & \hskip -1cm\\
 \hskip -0.05cm p_2\circ (((m\circ (1\times g))\times 1)\circ (1\times <1,1>))=\, !  & \hskip -1cm\endcases $ 
\vskip 0.0cm $\cases  \hskip -0.1cm*=p_1\circ ((m\circ (1\times g)\, \circ <1,!>)\times 1)=m\circ (1\times g)\, \circ <1,!>\circ \, p_1=m\circ (1\times g)\circ 1_{G\times 1}=*  &  \hskip -1cm\\
  \hskip 11.5cm*=m\circ (1\times g)    &  \hskip -1cm  \endcases $ 
\vskip 0.1cm
\item{}Proof for left shift $L_g$ is similar.   \hfill $\square $
\enddemo

\vskip 0.2cm
{\bf Corollary.} If $\xymatrix{K \ar@{>->}[r] & G}$ is a subgroup (object) of $G$ then $\bold{C}(1,K)@>\sigma >>\bold{Aut}_{\bold{C}}(G):k\mapsto R_k$
is a (right) action $(G,K,\sigma )$ on $G$ (by right shifts from $K$). If quotient $Q(G,K,\sigma )\in Ob\, \bold{C}$ exists it is called
{\bf quotient object} $G/K$ (under right action of $K$).        \hfill $\square $

\vskip 0.3cm
\proclaim{\bf Proposition 8.1.2} Let $p:G\to G/K$ be a quotient map s.t. $\bold{C}(1,p):\bold{C}(1,G)\to \bold{C}(1,G/K)$ is surjective. Then
\item{$\bullet $} $L_g:G\to G$ induces a $\bold{Set}$-map $\bar L'_g:\bold{C}(1,G/K)@>\sim >>\bold{C}(1,G/K)$  
\vskip -0.5cm
\item{$\bullet $} If $G/K$ is universal then $\exists \, \bar L_g:G/K@>\sim >> G/K$ in $\bold{C}$ s.t. $\bold{C}(1,\bar L_g)=\bar L'_g$
and $\vcenter{\xymatrix{G \ar[r]^-{L_g}_-{\sim } \ar[d]_-{p} &  G  \ar[d]^-{p} \\
G/K  \ar[r]_-{\bar L_g}^-{\sim } &   G/K}}$
\endproclaim
\demo{Proof} 
\vskip -0.7cm
\item{$\bullet $} Claim 1: $\bar L'_g:px\mapsto pL_gx$ is well-defined and iso \hskip 0.5cm$\vcenter{\xymatrix{G \ar[rr]^-{L_g} \ar[dd]^-{p} & & G \ar[dd]^-{p} \\
   & 1 \ar[ul]_-{x} \ar[ur]^-{L_gx} \ar[dl]_-{px} \ar@{-->}[dr]^-{pL_gx} &  \\
C/K  &  & G/K}}$
\item{}Proof of Claim 1: If $px=px'$ then $\exists \, k\in \bold{C}(1,K)$ such that $x'=R_kx=x\cdot k$.
Then $L_gx'=L_g(x\cdot k)=g\cdot (x\cdot k)=(g\cdot x)\cdot k=(L_gx)\cdot k=R_k(L_gx)$, i.e. $pL_gx'=pL_gx$.
\item{}$\bar L'_{g^{-1}}\circ \bar L'_g(px)=\bar L'_{g^{-1}}(pL_gx)=p(L_{g^{-1}}L_gx)=px$.  \hfill $\square $
\vskip 0.1cm
\item{$\bullet $} Claim 2: $G\times G@>1\times p>>G\times G/K$ is a quotient map of $(G\times G,K,<1,\sigma >)\in Ob\, \bold{C}\text{-}K$, 
where $\bold{C}\text{-}K$ is a category of right actions of $\bold{C}(1,K)$ on objects of $\bold{C}$, \, $<1,\sigma >:\bold{C}(1,K)\to \bold{Aut}_{\bold{C}}(G\times G):k\mapsto 1\times R_k$.
\vskip 0.1cm
\item{}Proof of Claim 2 follows immediately from the definition of universal quotient.    \hfill   $\square $
\vskip 0.1cm
\item{}Claim 3: $\forall \, k:1\to G$ \ $m\circ (1\times R_k)=R_k\circ m$.
\vskip 0.1cm
\item{}Proof of Claim 3 follows from the diagram \vskip 0.0cm
$\vcenter{\xymatrix{G\times G \ar@{..>}[rr]^(0.55){\sim }_(0.55){1\times p_1^{-1}} \ar@/^1.5pc/[rrrrr]^-{1\times R_k} \ar@{..>}[d]_-{m} & & G\times G\times 1 \ar[rr]_-{1\times 1\times k} \ar@{..>}[d]_-{m\times 1} & & G\times G\times G \ar[r]_-{1\times m} \ar[d]_-{m\times 1} & G\times G \ar[r]^-{m} & G \\
G \ar@{..>}[rr]^(0.55){\sim }_(0.55){p_1^{-1}} \ar@/_3.5pc/[urrrrrr]_-{R_k} & & G\times 1 \ar[rr]_-{1\times k} & & G\times G \ar[urr]_-{m} & & }}$
\item{}The dotted paths $p_1^{-1}\circ m, \, (m\times 1)\circ (1\times p_1^{-1}):G\times G\to G\times 1$ are equal since their composites with projections $p_1:G\times G\to G, \ p_2=\, !:G\times G\to 1$ are equal: 
$\cases p_1\circ (p_1^{-1}\circ m)=m & \hskip -0.2cm\\
p_2\circ (p_1^{-1}\circ m)=\, !	& \hskip -0.2cm\endcases $ and 
\item{}$\cases p_1\circ (m\times 1)\circ (1\times p_1^{-1})=m\times p_{G\times G}\circ (1\times p_1^{-1})=m\circ 1_{G\times G}=m & \hskip -0.2cm \\
p_2\circ (m\times 1)\circ (1\times p_1^{-1})=\, ! & \hskip -0.2cm\endcases $ 
\vskip 0.1cm
\item{}[where $p_{G\times G}:G\times G\times 1\to G\times G$ is the projection, $p_{G\times G}\circ (1\times p_1^{-1})=1_{G\times G}$ since  
$p_{G\times G}\circ (1\times p_1^{-1})\, \, \circ <x,y>=p_{G\times G}\, \, \circ <x,y,!>=<x,y>$]. \hfill $\square $
\vskip 0.0cm
\item{}Claim 4: $\exists \, !$ map $\bar m:G\times G/K\to G/K$ such that \hskip 0.35cm$\vcenter{\xymatrix{G\times G  \ar[r]^-{m} \ar[d]_-{1\times p} &  G  \ar[d]^-{p} \\
G\times G/K  \ar@{-->}[r]_-{\exists ! \bar m}  &    G/K }}$
\item{}Proof of Claim 4: Take $<f,g>\sim _{K}<f',g'>:X\to G\times G$ (where $G\times G$ is with the action $\bold{C}(1,K)\ni k\mapsto 1\times R_k\in \bold{Aut}_{\bold{C}}(G\times G)$)
then $<f',g'>=(1\times R_k)\, \circ <f,g>$ for some $k\in \bold{C}(1,K)$. 
\item{}$p\circ m\, \circ <f',g'>=p\circ m\circ (1\times R_k)\, \circ <f,g>=p\circ R_k\circ m\, \circ <f,g>=p\circ m\, \circ <f,g>$ 
[$p\circ R_k=p$ by definition of quotient $p$]. So, $p\circ m$ coequalizes $\sim _K$-equivalent arrows to $G\times G$. Therefore, 
$\exists ! \, \bar m:G\times G/K\to G/K$ filling out the above diagram.  \hfill  $\square $
\vskip 0.1cm
Now, define $\bar L_g:G/K\to G/K$ \, as \, $G/K@>\sim >>1\times G/K@>g\times 1>>G\times G/K@>\bar m>>G/K$. It works since
\item{}\hskip 1.5cm$\xymatrix{G \ar[rr]^(0.55){\sim } \ar@/^1.5pc/[rrrrrr]^-{L_g} \ar[d]_-{p} && 1\times G \ar[rr]^-{g\times 1} \ar[d]_-{1\times p} && G\times G \ar[rr]^(0.45){m} \ar[d]^-{1\times p} && G \ar[d]^-{p} \\
G/K \ar[rr]_(0.55){\sim } \ar@/_1.5pc/[rrrrrr]_-{\bar L_g} && 1\times G/K \ar[rr]_-{g\times 1} && G\times G/K \ar[rr]_(0.45){\bar m} && G/K }$
\item{}and \, $\bold{C}(1,\bar L_g)(px)=\bar L_g\circ p\circ x=p\circ L_g\circ x=pL_gx=\bar L'_g(px)$, i.e. $\bold{C}(1,\bar L_g)=\bar L'_g$.   \hfill $\square $
\enddemo

\vskip 0.2cm
\subhead 8.2. $\bold{C}$-group actions \endsubhead

\vskip 0.2cm
{\bf Definition 8.2.1.} 
\item{$\bullet $} Let $G$ be a group object in $\bold{C}$, $X\in Ob\, \bold{C}$, then $\bold{C}$-map $\rho :G\times X\to X$  
is a (left) {\bf group action} on $X$ if \hskip 1cm$\vcenter{\xymatrix{G\times G\times X \ar[r]^-{1\times \rho } \ar[d]_-{m\times 1} & G\times X \ar[d]^-{\rho } \\
G\times X \ar[r]_-{\rho } & X}}$ \hskip 1cm
$\vcenter{\xymatrix{1\times X \ar[r]^-{e\times 1} \ar[dr]_-{p_2} & G\times X \ar[d]^-{\rho } \\
     &   X }}$
\item{$\bullet $} {\bf Left shift} $L^X_g:X\to X$ (by $g\in \bold{C}(1,G)$) is the composite \hskip 0.25cm$\vcenter{\xymatrix{X \ar[r]_(0.5){\sim }^(0.53){p_2^{-1}} \ar[drr]_-{L^X_g} & 1\times X \ar[r]^-{g\times 1} & G\times X \ar[d]^-{\rho } \\
 & & X }}$
\item{$\bullet $} If $\xymatrix{K \ar@{>->}[r]^-{i_1} & G}$ is a subgroup of $G$, $\xymatrix{Y \ar@{>->}[r]^-{i_2} & X}$ is a subobject of $X$ 
then $K$ {\bf stabilizes} $Y$ if $\exists f:K\times Y\to Y$ such that \hskip 0.8cm$\vcenter{\xymatrix{G\times X \ar[r]^-{\rho } & X \\
K\times Y \ar[r]_-{f} \ar@{>->}[u]^-{i_1\times i_2} & Y \ar@{>->}[u]_-{i_2} }}$    \hfill $\square $

\vskip 0.1cm
\proclaim{\bf Lemma 8.2.1} Let $\xymatrix{Y \ar@{>->}[r]^-{i} & X}$ be a subobject of object $X$ with $G$-action $\rho :G\times X\to X$.
Assignment $\bold{Stab}_Y:Ob\, \bold{C}\to Ob\, \bold{Set}:Z\mapsto \bold{Stab}_Y(Z)\subset \bold{C}(Z,G)$ such that 
$(x:Z\to G)\in \bold{Stab}_Y(Z)$ \, iff \, $\exists \, \rho _x:Z\times Y\to Y$ such that \hskip 0.3cm$\vcenter{\xymatrix{G\times X \ar[r]^-{\rho } & X \\
G\times Y \ar@{>->}[u]^-{1\times i} & Y \ar@{>->}[u]_-{i} \\
Z\times Y \ar[u]^-{x\times 1} \ar@{-->}[ur]_-{\exists \rho _x} & }}$\hskip 0.3cm
is functorial (hom-subfunctor).
\endproclaim 
\demo{Proof} For $(f:W\to Z)\in Ar\, \bold{C}$ define $\bold{Stab}_Y(f):\bold{Stab}_Y(Z)\to \bold{Stab}_Y(W):x\mapsto x\circ f$ 
(as precomposite with $f$). This is correct since if $x\in \bold{Stab}_Y(Z)$ then $x\circ f\in \bold{Stab}_Y(W)$ 
which can be seen from the diagram \hskip 0.3cm$\vcenter{\xymatrix{G\times X \ar[r]^-{\rho } & X \\
Z\times Y \ar[r]^-{\rho _x} \ar[u]^-{x\times i} & Y \ar@{>->}[u]_-{i} \\
W\times Y \ar[u]^-{f\times 1} \ar@{-->}[ur]_-{\rho _{x\circ f}} & }}$\hskip 0.3cm
Functorial properties of $\bold{Stab}_Y$ are obvious. So, $\exists $ functor $\bold{Stab}_Y:\bold{C}^{op}\to \bold{Set}$ \, and \,  
$\bold{Stab}_Y\hookrightarrow \bold{C}(-,G)$ \, is a hom-subfunctor.    \hfill $\square $
\enddemo 

\vskip 0.2cm
{\bf Definition 8.2.2.} If $\bold{Stab}_Y:\bold{C}^{op}\to \bold{Set}$ is representable then denote its representing object by
\, $\text{\tt Stab}_Y\in Ob\, \bold{C}$ \, and call it a {\bf stabilizer} of $\xymatrix{Y \ar@{>->}[r] & X}$ (for group $G$ acting on $X$). \hfill $\square $

\vskip 0.2cm
\proclaim{\bf Proposition 8.2.1} Let \, $\bold{Stab}_Y:\bold{C}^{op}\to \bold{Set}$ be represented by \, $\text{\tt Stab}_Y\in Ob\, \bold{C}$. Then
\item{$\bullet $} $\xymatrix{\text{\tt Stab}_Y \ar@{>->}[r]^-{j} & G}$ is a subobject of group $G$ (but not necessarily a group object itself),
\item{$\bullet $} $j$ is the universal element of functor $\bold{Stab}_Y$,
\item{$\bullet $} each element in \, $\bold{Stab}_Y(Z)$ \, has form \, $j\circ x$ \, for a unique \, $x:Z\to \text{\tt Stab}_Y$, 
and all elements of this form ($\forall x:Z\to \text{\tt Stab}_Y$) are in $\bold{Stab}_Y(Z)$ [in other words, $(z\in _ZG)\, \& \, (z\in \bold{Stab}_Y(Z))\Leftrightarrow (z\in _Z\text{\tt Stab}_Y)$],
\item{$\bullet $} $\exists !\, \rho _j:\text{\tt Stab}_Y\times Y\to Y$ such that $\vcenter{\xymatrix{G\times X \ar[r]^-{\rho } & X \\
\text{\tt Stab}_Y\times Y \ar[u]^-{j\times i} \ar[r]_-{\rho _j} & Y \ar@{>->}[u]_-{i} }}$
and $\rho _j$ is universal among arrows with similar property, i.e. $\forall \rho _{x'}:Z\times Y\to Y$ such that 
$\vcenter{\xymatrix{G\times X \ar[r]^-{\rho } & X \\
Z\times Y \ar[u]^-{x'\times i} \ar[r]_-{\rho _{x'}} & Y \ar@{>->}[u]_-{i} }}$
$\exists !\, x:Z\to \text{\tt Stab}_Y$ such that $x'=j\circ x$ and
$\vcenter{\xymatrix{\text{\tt Stab}_Y\times Y \ar[r]^-{\rho _j} & Y \\
Z\times Y \ar[u]^-{x\times 1} \ar[ur]_-{\rho _{x'}} & }}$
\endproclaim 
\demo{Proof} First three points follow from Yoneda Lemma ($j$ is the universal element of representation $\bold{C}(-, \text{\tt Stab}_Y)@>\sim >>\bold{Stab}_Y$ corresponding 
(under Yoneda embedding) to monic (natural transformation) $\xymatrix{\bold{C}(-,\text{\tt Stab}_Y) \ar[r]^-{\sim } \ar@{>->}@/_1.0pc/[rr] & \bold{Stab}_Y \ar@{^{(}->}[r] & \bold{C}(-,G)}$).
Fourth point follows from the definition and above properties of functor $\bold{Stab}_Y$ and that $\xymatrix{Y \ar@{>->}[r]^-{i} & X }$ is monic.   \hfill $\square $ 
\enddemo 

\vskip 0.2cm
\proclaim{\bf Lemma 8.2.2} Subobject $\xymatrix{H \ar@{>->}[r] & G}$ of a group object $G\in Ob\, \bold{C}$ is itself a group object iff
$\forall Z\in Ob\, \bold{C}$ there are induced group operations in hom-set $\bold{C}(Z,H)$ in the following way
$\xymatrix{G\times G \ar[r]^-{m} & G \\
H\times H \ar@{>->}[u] & H \ar@{>->}[u] \\
Z \ar[u]^-{<x,y>} \ar@{-->}[ur]_-{\exists \, m(<x,y>)} & }$\hskip 1.5cm
$\xymatrix{G \ar[r]^-{inv} & G \\
H \ar@{>->}[u] & H \ar@{>->}[u] \\
Z \ar[u]^-{x} \ar@{-->}[ur]_-{\exists \, inv(x)} & }$\hskip 1.5cm
$\xymatrix{1 \ar[r]^-{e} & G \\
Z \ar[u]^-{!} \ar@{-->}[r]_-{\exists \, e(!)} & H \ar@{>->}[u] }$
\endproclaim 
\demo\nofrills{Proof \ } is obvious in both directions.    \hfill $\square $
\enddemo 

\vskip 0.2cm
\proclaim{\bf Proposition 8.2.2} 
\item{$\bullet $} $\xymatrix{\text{\tt Stab}_Y \ar@{>->}[r] & G}$ is always a submonoid of group object $G\in Ob\, \bold{C}$.
\item{$\bullet $} $\xymatrix{\text{\tt Stab}_Y \ar@{>->}[r] & G}$ is a subgroup of group object $G\in Ob\, \bold{C}$ if 
$\forall Z\in Ob\, \bold{C} \ \forall x\in \bold{Stab}_Y(Z)$ the corresponding map \, $\rho _x:Z\times Y\to Y$ as in the diagram 
\hskip 0.2cm$\vcenter{\xymatrix{G\times X \ar[r]^-{\rho } & X \\
Z\times Y \ar[u]^-{x\times i} \ar@{-->}[r]_-{\rho _x} & Y \ar@{>->}[u]_-{i} }}$\hskip 0.2cm
is surjective in the second argument, i.e. \, $\forall t:T\to Z$ \, the map \, $\bold{C}(T,Y)\ni s\mapsto \rho _x\circ <t,s>\, \in \bold{C}(T,Y)$
is surjective [it holds in classical case in $\bold{Set}$].
\endproclaim 
\demo{Proof} 
\vskip 0.1cm
\item{$\bullet $} \hskip -0.3cm$
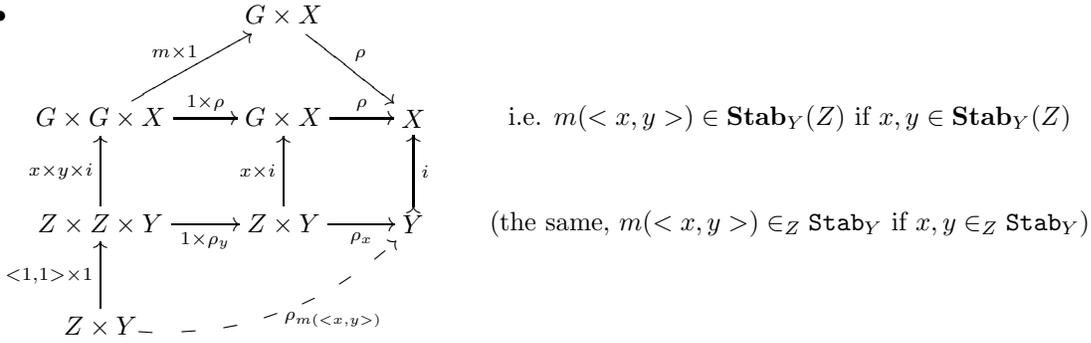
\xymatrix{ & G\times X \ar[dr]^-{\rho } & & \\
G\times G\times X \ar[ur]^-{m\times 1} \ar[r]^-{1\times \rho } & G\times X \ar[r]^-{\rho } & X & {\hskip -0.2cm\text{i.e. }{m(<x,y>)\in \bold{Stab}_Y(Z)} \text{ if } {x,y\in \bold{Stab}_Y(Z)}} \\
Z\times Z\times Y \ar[u]^-{x\times y\times i} \ar[r]_-{1\times \rho _y} & Z\times Y \ar[u]^-{x\times i} \ar[r]_-{\rho _x} & Y \ar@{>->}[u]_-{i} & {\hskip -0.2cm(\text{the same, } {m(<x,y>)\in _Z\text{\tt Stab}_Y} \text{ if } {x,y\in _Z\text{\tt Stab}_Y})} \\
Z\times Y \ar[u]^-{<1,1>\times 1} \ar@/_1.35pc/@{-->}[urr]_(0.6){\ \rho _{m(<x,y>)}} & & & }$
\vskip 0.2cm
\hskip 0.35cm$
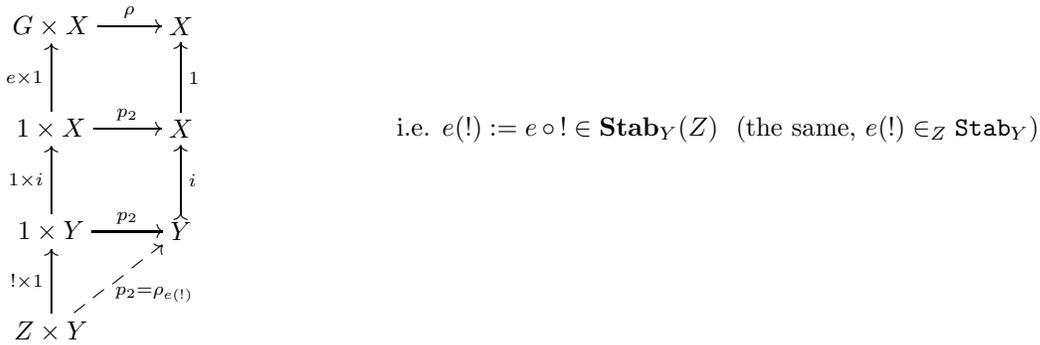
\xymatrix{G\times X \ar[r]^-{\rho } & X & & \\
1\times X \ar[u]^-{e\times 1} \ar[r]^-{p_2} & X \ar[u]_-{1} & & \hskip 0.6cm\text{i.e. } e(!):=e\, \circ \, !\in \bold{Stab}_Y(Z) \ \ (\text{the same, } e(!)\in _{Z}\text{\tt Stab}_Y) \\
1\times Y \ar[u]^-{1\times i} \ar[r]^-{p_2} & Y \ar@{>->}[u]_-{i} & & \\
Z\times Y \ar[u]^-{!\times 1} \ar@{-->}[ur]_-{p_2=\rho _{e(!)}} & & & }$
\vskip 0.2cm
\item{$\bullet $} In general, for \, $x\in _Z\text{Stab}_Y$ \, $inv(x)\in _ZG$, \, but \, $inv(x)\notin _Z\text{\tt Stab}_Y$.
\item{} \hskip 0.5cm$
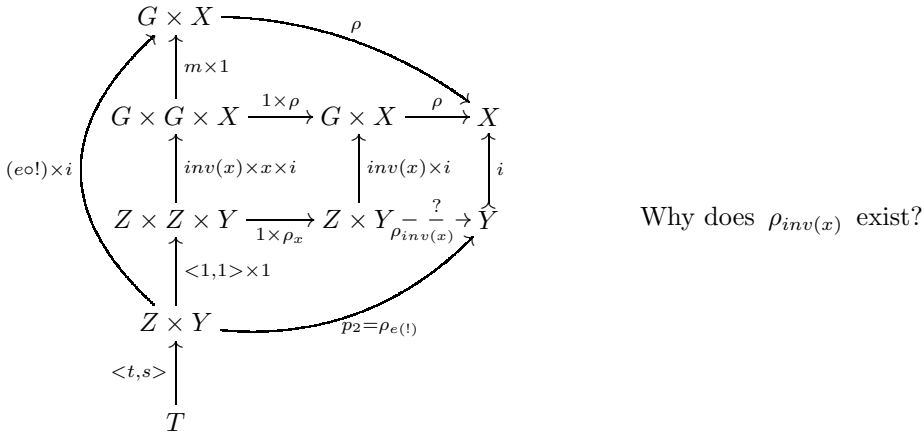
\xymatrix{G\times X \ar@/^1pc/[drr]^-{\rho } & & & \\
G\times G\times X \ar[u]_-{m\times 1} \ar[r]^-{1\times \rho } & G\times X \ar[r]^-{\rho } & X & \\
Z\times Z\times Y \ar[u]_-{inv(x)\times x\times i} \ar[r]_-{1\times \rho _x} & Z\times Y \ar[u]_-{inv(x)\times i} \ar@{-->}[r]^-{?}_(0.495){\rho _{inv(x)}} & Y \ar@{>->}[u]_-{i} & \hskip 0.8cm\text{Why does \, } \rho _{inv(x)} \text{ \, exist?}\\
Z\times Y \ar[u]_-{<1,1>\times 1} \ar@/_1.5pc/[urr]_-{\ p_2=\rho _{e(!)}} \ar@/^3pc/[uuu]^-{(e\circ !)\times i} & & & \\
T \ar[u]^-{<t,s>} & & & }$

\vskip 0.2cm
\proclaim{\bf Lemma 8.2.3} If a square \hskip 0.3cm$\vcenter{\xymatrix{Hom(-,A) \ar@{->>}[r]^-{f} \ar[d]_-{h} & Hom(-,B) \ar[d]^-{g} \ar@{-->}[dl]^-{\exists \, d} \\
Hom(-,C) \ar@{>->}[r]_-{k} & Hom(-,D)}}$\hskip 0.3cm of natural transformations of representables commutes, where \, $f$ \, and \, $k$ \, are respectively 
componentwise surjective and componentwise injective maps, then \, $\exists $ \, a (unique) diagonal \, $d$ \, keeping the diagram commutative. 
\endproclaim 
\demo{Proof of Lemma} \hskip 2cm$\xymatrix{Hom(X,A) \ar@{->>}[rr]^-{f_X} \ar[dd]_-{h_X} & & Hom(X,B) \ar[dd]^-{g_X} \ar[dl]_-{g'_X} \\
     & \hskip -0.15cmIm(k) \ar@{^{(}->}[dr] &     \\
Hom(X,C) \ar[ur]^-{k'_X}_-{\sim } \ar@{>->}[rr]_-{k_X} & & Hom(X,D) }$
\item{}$k_X,\, g_X$ factor through \, $Im(k)$ \, since \, $k_X$ \, is injective and \, $f_X$ \, is surjective. Define diagonal \, $d_X:=(k'_X)^{-1}\circ g'_X$. 
Arrows \, $d_X$ ($X$ is a parameter) form natural transformation which can be seen from the diagram 
\hskip 2.5cm$\xymatrix{ & & \cdot \ar[r] \ar[d] \ar[dll] & \cdot \ar[d] \ar@{..>}[dl]^-{\hskip -0.15cmd_X} \ar@{..>}[dll] \\
\cdot \ar[r] \ar[d] & \cdot \ar[d] \ar@{..>}[dl]_-{d_{X'}\hskip -0.175cm} & \cdot \ar[r] \ar@{..>}[dll] & \cdot \ar[dll] \\
\cdot \ar@{>->}[r]_-{monic} & \cdot & & }$      \vskip -2cm   \hfill  $\square $
\enddemo 
\vskip 1.5cm
So, apply the above lemma to the square \hskip 0.35cm
$\vcenter{\xymatrix{\bold{C}(-,Z\times Y) \ar@{->>}[rr]^-{(1\times \rho _x)\circ (<1,1>\times 1)} \ar[d]_-{p_2} & & \bold{C}(-,Z\times Y) \ar[d]^-{\rho \circ (inv(x)\times i)} \ar@{-->}[dll]^-{\ \ \ \, \exists \, \rho _{inv(x)}} \\
\bold{C}(-,Y) \ar@{>->}[rr]_-{i} & & \bold{C}(-,X) }}$
\item{}(The top arrow is componentwise surjective since $(1\times \rho _x)\circ (<1,1>\times 1)\, \circ <t,s>=(1\times \rho _x)\, \circ <~t,t,s>=<t,\rho _x\, \circ <t,s>>$, 
and, so that, $\forall <m,l>:T\to Z\times Y$ \, $\exists $ \, its preimage \, $<t,s>:T\to Z\times Y$ \, with \, $t=m$ \, and \, $s$ \, is a solution of the equation \, $\rho _x\, \circ <m,s>=l$ \, [which exists because \, $\rho _x$ \, is surjective in the second argument].
The bottom arrow is componentwise injective since \, $i$ \, is monic.)
\item{}Consequently, $\exists (!) \, \rho _{inv(x)}:Z\times Y\to Y$ such that $\vcenter{\xymatrix{G\times X \ar[r]^-{\rho } & X \\
Z\times Y \ar[u]^-{inv(x)\times i} \ar[r]_-{\rho _{inv(x)}} & Y \ar@{>->}[u]_-{i} }}$ i.e. 
if \, $x\in _Z\text{\tt Stab}_Y$ then \, $inv(x)\in _Z\text{\tt Stab}_Y$.    \hfill   $\square $
\enddemo 

\vskip 0.2cm
\proclaim{\bf Lemma 8.2.4} If $L_g:X@>\sim >>X$ is a left shift, and $\xymatrix{Y,Z \ar@{>->}[r] & X}$ are subobjects of $X$ 
such that an induced isomorphism $\overline L_g:Y@>\sim >>Z$ exists, i.e. the diagram $\vcenter{\xymatrix{X \ar[r]^-{L_g} & X \\
Y \ar@{>->}[u]^-{i_Y} \ar@{-->}[r]_-{\exists \, \overline L_g}^-{\sim } & Z \ar@{>->}[u]_-{i_Z} }}$ commutes, then $\exists $ an induced map (iso)
$\overline {L_g\circ R_{g^{-1}}}:\text{\tt Stab}_Y@>\sim >>\text{\tt Stab}_Z$, corresponding to $\overline L_g:Y@>\sim >>Z$, 
such that the following diagram commutes\hskip 0.39cm
$\vcenter{\xymatrix{G\times X \ar[rrr]^-{\rho } & & & X & \\
 & G\times X \ar[rrr]^-{\rho } \ar[ul]_(0.35){\ \ \ (L_g\circ R_{g^{-1}})\times L_g} & & & X \ar[ul]_-{\hskip -0.1cmL_g} \\
\text{\tt Stab}_Z\times Z \ar@{>->}[uu]^-{j_Z\times i_Z} \ar@{..>}[rrr]^(0.65){\rho _{j_Z}} & & & Z \ar@{>.>}[uu]_(0.35){i_Z} & \\
 & \text{\tt Stab}_Y\times Y \ar[ul]^-{\overline {L_g\circ R_g^{-1}}\times \overline L_g\hskip 0.3cm} \ar@{>->}[uu]^(0.725){j_Y\times i_Y\hskip -0.05cm} \ar[rrr]_-{\rho _{j_Y}} & & & Y \ar@{..>}[ul]_-{\hskip -0.1cm\overline L_g}^(0.6){\sim \hskip -0.05cm} \ar@{>->}[uu]_-{i_Y} }}$
\endproclaim
\demo{Proof} The only difficulty is that the left side square commutes, namely, $\vcenter{\xymatrix{G & G \ar[l]_-{L_g\circ R_{g^{-1}}} \\
\text{\tt Stab}_Z \ar@{>->}[u]^-{j_Z} & \text{\tt Stab}_Y \ar@{>->}[u]_-{j_Y} \ar@{-->}[l]^-{\underset {\overline {L_g\circ R_{g^{-1}}}}\to {}}_-{?} \hskip -0.2cm}}$
Sufficient to show that $(L_g\circ R_{g^{-1}})\circ j_Y\in _{\underset {}\to {\text{\tt Stab}_Y}}\text{\tt Stab}_Z$, i.e. that 
$(L_g\circ R_{g^{-1}})\circ j_Y\in \bold{Stab}_Z(\text{\tt Stab}_Y)$, or that $\exists $\, $\rho ':\text{\tt Stab}_Y\times Z\to Z$ such that
\hskip 0.8cm$\vcenter{\xymatrix{G\times X \ar[r]^-{\rho } & X \\
\text{\tt Stab}_Y\times Z \ar@{-->}[r]_-{\rho '} \ar[u]^-{((L_g\circ R_{g^{-1}})\circ j_Y)\times i_Z} & Z \ar@{>->}[u]_-{i_Z} }}$
\hskip 0.8cmIndeed, $((L_g\circ R_{g^{-1}})\circ j_Y)\times (L_g\circ i_Y)=((L_g\circ R_{g^{-1}})\circ j_Y)\times (i_Z\circ \overline L_g)=(((L_g\circ R_{g^{-1}})\circ j_Y)\times i_Z)\circ (1\times \overline L_g)$,
then $\rho \circ (((L_g\circ R_{g^{-1}})\circ j_Y)\times i_Z)\circ (1\times \overline L_g)=i_Z\circ \overline L_g\circ \rho _{j_Y}$, 
and so, $\rho \circ (((L_g\circ R_{g^{-1}})\circ j_Y)\times i_Z)=i_Z\circ \overline L_g\circ \rho _{j_Y}\circ (1\times \overline L_g)^{-1}$, i.e.
$\rho ':=\overline L_g\circ \rho _{j_Y}\circ (1\times \overline L_g)^{-1}$. 
Therefore, $\forall x\in _T\text{\tt Stab}_Y$ \, $L_g\circ R_{g^{-1}}(x)\in _T\text{\tt Stab}_Z$, i.e. 
$\exists $ the induced map $\overline {L_g\circ R_{g^{-1}}}:\text{\tt Stab}_Y\to \text{\tt Stab}_Z$.      \hfill    $\square $
\enddemo

\vskip 0.2cm
\proclaim{\bf Proposition 8.2.3} 
\item{$\bullet $} Two objects $(G,\bold{C}(1,\text{\tt Stab}_Y),\sigma _1)$ and $(G,\bold{C}(1,\text{\tt Stab}_Z),\sigma _2)$
from $\bold{C}\text{-}\bold{Grp}$ (a category of right group actions on objects from $\bold{C}$) are (equivariantly) isomorphic
if $\exists g\in \bold{C}(1,G)$ and an induced isomorphism (as in lemma 2.7.2.4) $\overline L_g:Y@>\sim >>Z$. The required isomorphism
has form $\xymatrix{(G,\bold{C}(1,\text{\tt Stab}_Y),\sigma _1) \ar[rrr]_-{\sim }^-{(L_g\circ R_{g^{-1}},\, \overline {L_g\circ R_{g^{-1}}}\circ -)} & & & (G,\bold{C}(1,\text{\tt Stab}_Z),\sigma _2) }$\hskip -0.1cm.
\vskip 0.05cm
\item{$\bullet $} $G/\text{\tt Stab}_Y\simeq G/\text{\tt Stab}_Z$ (if these quotients exist).
\endproclaim 
\demo{Proof} 
\item{$\bullet $} It is necessary to prove that $\forall g:1\to G$ and $k:1\to \text{\tt Stab}_Y$ \, $L_g\circ R_{g^{-1}}\circ R_k=R_{g\cdot k\cdot g^{-1}}\circ L_g\circ R_{g^{-1}}$. 
It follows from two facts $R_{g\cdot k\cdot g^{-1}}=R_{g^{-1}}\circ R_k\circ R_g$ (antihomomorphism) and commutativity of left and right 
shifts $L_{g_1}\circ R_{g_2}=R_{g_2}\circ L_{g_1}$ [the last fact follows from associativity axiom 
$\forall <t,r,s>:T\to G\times G\times G$ \, $m(t,m(r,s))=m(m(t,r),s)$, and so, $L_{g_1}\circ R_{g_2}\circ t=m(g_1\circ !,m(t,g_2\circ !))=m(m(g_1\circ !,t),g_2\circ !)=R_{g_2}\circ L_{g_1}\circ t$]. 
\item{$\bullet $} Isomorphic objects in $\bold{C}\text{-}\bold{Grp}_Q$ have isomorphic quotients in $\bold{C}$ since 
$Q:\bold{C}\text{-}\bold{Grp}_Q\to \bold{C}$ is a functor. So, $Q(G,\bold{C}(1,\text{\tt Stab}_Y),\sigma _1)\simeq Q(G,\bold{C}(1,\text{\tt Stab}_Z),\sigma _2)$.    \hfill   $\square $
\enddemo 

\vskip 0.15cm
{\bf Definition 8.2.3.} An object $X\in Ob\, \bold{C}$ with a group action $\rho :G\times X\to X$ such that $\forall x:1\to X$
both $\text{\tt Stab}_x$ and $G/\text{\tt Stab}_x$ exist, and $G/\text{\tt Stab}_x$ is universal, is called {\bf homogenious} 
if $\exists $ an isomorphism $f:G/\text{\tt Stab}_x@>\sim >>X$ such that
$\vcenter{\xymatrix{G\times (G/\text{\tt Stab}_x) \ar[r]^-{\rho '} \ar[d]^-{\sim }_-{1\times f} & G/\text{\tt Stab}_x \ar[d]_-{\sim }^-{f} \\
G\times X \ar[r]_-{\rho } & X }}$ (for an $x:1\to X$)
where $\rho '$ is defined from $\vcenter{\xymatrix{G\times G \ar[r]^-{m} \ar[d]_-{1\times p} & G \ar[d]^-{p} \\
G\times (G/\text{\tt Stab}_x) \ar@{-->}[r]_-{\exists !\, \rho '} & G/\text{\tt Stab}_x }}$
($1\times p$ and $p$ are quotient maps).      \hfill   $\square $

\vskip 0.2cm
\proclaim{\bf Proposition 8.2.4} If $X$ is a homogenious object (with $G$-action $\rho :G\times X\to X$) and
$\bold{C}(1,p):\bold{C}(1,G)\to \bold{C}(1,G/\text{\tt Stab}_x)$ is surjective, where $\xymatrix{G \ar@{->>}[r]^-{p} & G/\text{\tt Stab}_x }$
is a quotient map, then 
\item{$\bullet $} $\bold{C}(1,G)$ acts transitively on $\bold{C}(1,X)$, i.e. $\forall x,y:1\to X$ $\exists g:1\to G$ such that 
$y=L_g\circ x$,
\item{$\bullet $} definition of homogenious object $X@>\sim >f^{-1}>G/\text{\tt Stab}_x$ does not depend on the choice of \, $x:1\to X$.
\endproclaim 
\demo{Proof}
\item{$\bullet $} $\forall a',b':1\to G/\text{\tt Stab}_x$ $\exists a,b:1\to G$ s.t. $pa=a', \, pb=b'$, and $\exists g:1\to G$
s.t. $b=L_ga$. By proposition 8.1.2, \, $\bar L'_g(a')=\bar L'_g(pa)=pL_g(a)=pb=b'$ (where $\bar L'_g$ is the induced left shift
on $\bold{C}(1,G/\text{\tt Stab}_x)$). So, $\bold{C}(1,G)$ acts transitively on set $\bold{C}(1,G/\text{\tt Stab}_x)$, and consequently, on $\bold{C}(1,X)$.
\item{$\bullet $}  
Regard the diagram \hskip 0.3cm
$\vcenter{\xymatrix{ & & G\times G \ar[r]^-{m} \ar@{..>}[d]_(0.35){1\times p_y\hskip -0.1cm} & G \ar[d]^-{p_y} \\
G\times G \ar[r]^-{m} \ar[d]_-{1\times p_x} \ar[urr]^-{(L_g\circ R_{g^{-1}})\times (L_g\circ R_{g^{-1}})\hskip 0.75cm} & G \ar[d]^-{p_x} \ar[urr]_(0.7){L_g\circ R_{g^{-1}}} & G\times (G/\text{\tt Stab}_y) \ar@{..>}[r]_-{\rho ''} & G/\text{\tt Stab}_y \\
G\times (G/\text{\tt Stab}_x) \ar[r]_-{\rho '}  \ar[d]_-{1\times f}^-{\sim } & G/\text{\tt Stab}_x  \ar[d]^(0.3){\hskip -0.05cmf}_(0.35){\sim } & & \\
G\times X \ar[r]_-{\rho } \ar@/_1.75pc/@{..>}[uurr]^(0.85){(L_g\circ R_{g^{-1}})\times \alpha \hskip -0.18cm}_(0.9){\hskip -0.1cm\sim } & X \ar@{-->}[uurr]_-{\exists !\, \alpha }^-{\sim } & & }}$
\item{}[\, $\alpha $ exists as a mediating arrow because $f\circ p_x$ is a quotient map of $(G,\bold{C}(1,\text{\tt Stab}_x),\sigma _1)$, and 
$p_y\circ (L_g\circ R_{g^{-1}})$ coequalizes $\sim _{\sigma _1}$-equivalent arrows ($L_g\circ R_{g^{-1}}$ is equivariant, 
and $p_y$ is a quotient map of $(G,\bold{C}(1,\text{\tt Stab}_y),\sigma _2)$), essentially, $\alpha =Q(L_g\circ R_{g^{-1}})$. 
The bottom square commutes because $(1\times f)\circ (1\times p_x)$ is a quotient map, and so, epi\, ]
\item{}Therefore, \hskip 0.3cm
$\vcenter{\xymatrix{G\times (G/\text{\tt Stab}_y) \ar[r]^-{\rho ''} & G/\text{\tt Stab}_y \\
G\times X \ar[r]^-{\rho } \ar[u]^-{(L_g\circ R_{g^{-1}})\times \alpha }_-{\sim } & X \ar[u]_-{\alpha }^-{\sim } \\
G\times X \ar[r]_-{\rho } \ar[u]^-{(L_{g^{-1}}\circ R_g)\times L_{g^{-1}}}_-{\sim } \ar@/^11pc/[uu]_-{\hskip 0.1cm1\times (\alpha \, \circ L_{g^{-1}})}^-{\sim } & X \ar[u]_-{L_{g^{-1}}}^-{\sim } }}$
\hskip 0.3cmi.e. $L_g\circ \alpha ^{-1}$ is the required isomorphism (by definition 8.2.3).    \hfill    $\square $
\enddemo

\head {\bf Bibliography} \endhead

\vskip 0.1cm
\refstyle{A}
\widestnumber\key{AAAAA}

\ref\key A-H-S
\by J. Adamek, H. Herrlich, G.E. Strecker
\book Abstract and Concrete Categories. The Joy of Cats
\yr online edition, 2004
\endref

\ref\key A-V-L
\by D.V. Alekseevskiy, A.M. Vinogradov, V.V. Lychagin
\book Main Ideas and Concepts of Differential Geometry
\yr 1988
\publ Moscow
\lang Russian
\endref

\ref\key Bi-Cr
\by R.L. Bishop, R.J. Crittenden
\book Geometry of Manifolds
\publ Academic Press, New York and London
\yr 1964
\endref

\ref\key Bor1
\by F. Borceux
\book Handbook of Categorical Algebra 1. Basic Category Theory
\yr 1994
\publ Cambridge University Press
\endref

\ref\key Bor2
\by F. Borceux
\book Handbook of Categorical Algebra 2. Categories and Structures
\yr 1994
\publ Cambridge University Press
\endref

\ref\key Bor3
\by F. Borceux
\book Handbook of Categorical Algebra 3. Categories of Sheaves
\yr 1994
\publ Cambridge University Press
\endref

\ref\key Bru
\by U. Bruzzo
\book Introduction to Algebraic Topology and Algebraic Geometry
\yr 2002
\publ International School for Advanced Studies, Trieste
\endref

\ref\key C-C-L
\by S.S. Chern, W.H.Chen, K.S. Lam
\book Lectures on Differential Geometry
\yr 2000
\publ World Scientific
\endref

\ref\key D-N-F
\by B.A. Dubrovin, S.P. Novikov, A.T. Fomenko
\book Modern Geometry
\publ Moscow
\yr 1979
\lang Russian
\endref

\ref\key Eng
\by R. Engelking
\book General Topology
\yr 1977
\endref

\ref\key ELOS
\by L.E. Evtushik, U.G. Lumiste, N.M. Ostianu, A.P. Shirokov
\book Differential Geometric Structures on Manifolds
\yr 1979
\publ Moscow
\lang Russian
\endref

\ref\key Hir
\by F. Hirzebruch
\book Topological Methods in Algebraic Geometry
\publ Springer-Verlag
\yr 1966
\endref

\ref\key Jac
\by B. Jacobs
\book Categorical Logic and Type Theory
\publ Elsevier, North-Holland
\yr 2001
\endref

\ref\key Kel
\by G.M. Kelly
\book Basic Concepts of Enriched Category Theory
\publ Reprints in TAC, No. 10
\yr 2005
\endref

\ref\key Kom
\by B.P. Komrakov
\book Structures on Manifolfs and Homogenious Spaces
\publ Minsk
\lang Russian
\yr 1978
\endref

\ref\key Mac
\by S. MacLane
\book Categories for the Working Mathematician
\publ Springer-Verlag
\yr 1971
\endref

\ref\key M-M
\by S. MacLane, I. Moerdijk
\book Sheaves in Geometry and Logic
\publ Springer-Verlag
\yr 1992
\endref

\ref\key May
\by J.P. May
\book A Concise Course in Algebraic Topology
\publ The University of Chicago Press
\yr 1999
\endref

\ref\key Moe
\by I. Moerdijk
\book Introduction to the Language of Stacks and Gerbes
\publ University of Utrecht, arXiv:math.AT/0212266v1
\yr 2002
\endref

\ref\key Nar
\by R. Narasimhan
\book Analysis on Real and Complex Manifolds
\publ North-Holland
\yr 1968
\endref

\ref\key Nes
\by J. Nestruev
\book Smooth manifolds and observables
\publ Moscow
\yr 2003
\lang Russian
\endref

\ref\key Str
\by T. Streicher
\book Fibred Categories a la Jean Benabou
\publ online lecture notes
\yr 1999
\endref

\ref\key Strt
\by R. Street
\paper Categorical and Combinatorial Aspects of Descent Theory
\publ talk at ICIAM
\yr 2003
\endref 

\ref\key Vis
\by A. Vistoli
\book Notes on Grothendieck Topologies, Fibered Categories and Descent Theory
\publ Bologna, Italy, online lecture notes
\endref

\enddocument